\documentclass{article}
\usepackage{amsmath,amsthm}
\usepackage{graphicx,epstopdf,epsfig,multirow,epic,bm}

\normalsize \rm
\begin{document}

\begin{center}
{\Large  \textbf { A Method for Geodesic Distance on Subdivision of Trees with Arbitrary Orders and Their Applications }}\\[12pt]
{\large Fei Ma$^{a,}$\footnote{~The author's E-mail: mafei123987@163.com. },\quad Ping Wang$^{b,c,d,}$\footnote{~The corresponding author's E-mail: pwang@pku.edu.cn.} \quad  and  \quad  Xudong  Luo$^{e,}$\footnote{~The author's E-mail: luoxudong117@163.com.} }\\[6pt]
{\footnotesize $^{a}$ School of Electronics Engineering and Computer Science, Peking University, Beijing 100871, China\\
$^{b}$ National Engineering Research Center for Software Engineering, Peking University, Beijing, China\\
$^{c}$ School of Software and Microelectronics, Peking University, Beijing  102600, China\\
$^{d}$ Key Laboratory of High Confidence Software Technologies (PKU), Ministry of Education, Beijing, China\\
$^{e}$ School of Mathematics and Statistics, Lanzhou University, 730000 Lanzhou, CHINA.}\\[12pt]
\end{center}

\begin{quote}
\textbf{Abstract:} Geodesic distance, sometimes called shortest path length, has proven useful in a great variety of applications, such as information retrieval on networks including treelike networked models. Here, our goal is to analytically determine the exact solutions to geodesic distances on two different families of growth trees which are recursively created upon an arbitrary tree $\mathcal{T}$ using two types of well-known operations, first-order subdivision and  ($1,m$)-star-fractal operation. Different from commonly-used methods, for instance, spectral techniques, for addressing such a problem on growth trees using a single edge as seed in the literature, we propose a novel method for deriving closed-form solutions on the presented trees completely. Meanwhile, our technique is more general and convenient to implement compared to those previous methods mainly because there are not complicated calculations needed.  In addition, the closed-form expression of mean first-passage time ($MFPT$) for random walk on each member in tree families is also readily obtained according to connection of our obtained results to effective resistance of corresponding electric networks. The results suggest that the two topological operations above are sharply different from each other due to $MFPT$ for random walks, and, however, have likely to show the similar performance, at least, on geodesic distance.


\end{quote}

\vskip 1cm

\section{Introduction}

Geodesic distance, conventionally called shortest path length in the language of graph theory, has proven useful in a great variety of areas, for example, discrete applied mathematics, theoretical computer science, biology science and so forth. While this is in essence not a new concept, its related researches still keep active in present science community and have been receiving more attention \cite{Petrushevski-2019,Yen-2013}. Included examples have information retrieval on internet \cite{Spyros-2017}, signal integrity in communication networks \cite{Petrushevski-2019}, disease spreading on relationship networks among individuals \cite{Jia-2018}, navigation in spatial networks \cite{Wei-2014}, to name but a few. Perhaps, one of the most important reasons for this is that it plays a helpful role for effectively measuring information spread and information retrieve on networks mentioned above. Specifically, for a pair of vertices $u$ and $v$ attempting to interact with one another on a network which is also abstractly thought of as a graph denoted by $\mathcal{G}(\mathcal{V},\mathcal{E})$, the much smaller the geodesic distance between them, the much less the time to need. Here, geodesic distance between vertices $u$ and $v$, defined as $d_{uv}$, is equal to the edge number of a shortest path connecting this pair of vertices. In addition, this concept is becoming increasingly popular in some newborn disciplines, for instance, complex network \cite{Brodka-2011}.

The recent two decades have seen a bloom of complex network study mainly because of its own right, which has helped us to better understand some synthetic and real-world complex systems. Such complex systems include World Wide Web (WWW), citation networks, metabolic networks, protein-protein interaction networks and predator-prey webs \cite{Michael-2019}-\cite{Ulrich-2019}. Particularly, the portion of complex network study is to focus mainly on investigating how dynamic and function taking place on complex networks of great interest are influenced by the topological structures of their underlying graphs. Among of which, one type of networked models have attracted significant attention according to their own specific topological structure like tree. One of such examples is citation networks of scientific papers \cite{Julia-2019} which is in general considered directed where the vertices represent documents and the directed edges represent citations between documents. In this respect, as one of the best studied models in various fields, traditional random graph can be used to generate networks of this type that are locally treelike, meaning that all local neighborhoods take the form of trees in the limit of large graph size. To be more general, a great deal of tree models have been directly proposed to account for some behaviors occurring in many complex networks \cite{MaF-2019}-\cite{Yuan-2010}. Among of them, a fraction of tree family become of great concern according to their own some intriguing topological structures including power-law distribution for vertex degree \cite{MaF-2019}, exponential vertex degree distribution \cite{Luo-2019}, fractal feature \cite{Peng-2014,Yuan-2010} which can turn out to be quite prevailing in man-made and natural complex networks.

As known, tree, the most fundamental and simplest connected graph, has been widely studied in the last \cite{Bondy-2008}-\cite{Zhang-2011}. Indeed, there are a surprising number of applications built upon structure of tree just because it is exactly solvable and owns  significant availability. Social networks, for instance, display hierarchical structure such that ones can naturally employ treelike models, the so-called dendrograms, as a graphical representation and summary of the structure of this type \cite{Aaron-2008}. Meanwhile, in information science, tree has been playing an important role in search engine at present. It can be worth mentioning that almost all data structures and a large fraction of algorithms suit for searching information are based on treelike models, such as binary tree \cite{Stefano-2018}. Among other things there are many other potential and invaluable applications on the basis of treelike models in coding theory including the best known Huffman tree coding, quedreecoding and so on \cite{Zhangj-2014}. The last but not least, some other interesting topics correlated with treelike models have received considerable attention, such as geodesic distance \cite{Deng-2019}, mean first-passage time for random walk \cite{Peng-2018}, fractal phenomena \cite{Nobutoshi-2019}, etc. We, in this paper, aim to study tree networks. To this end, we propose two different families of tree models, and analyze geodesic distance on them using a novel method developed later.

Roughly speaking, a central problem to answer is to explicitly capture solutions for geodesic distance on treelike models in view of theoretical flavor and practical applications. Some related works have been reported in the literature \cite{Deng-2019}. By far, the exact solutions for several types of tree models have been obtained by taking advantage of some typical tools in which the best used is Laplacian spectral and eigenvectors based on laplacian matrix of underlying structure \cite{Zhang-2011}. Nonetheless, we here do not adopt methods of this type but instead introduce some novel methods according to topological structure of tree models under consideration sufficiently. While the core of our methods is to built up a series of equations in an iterative manner akin to other enumeration methods \cite{Yuan-2010}, they are, in some extent, more convenient to manipulate, at least on tree models considered in this paper. More details about comparison between our work and early studies are deferred to show in Section 5. Our contributions are as follows.

(1) In order to obtain closed-form solution to geodesic distance on each member in tree families built later, we propose a novel computational method based on mapping which is more convenient than previous tools including spectral techniques.

(2) To highlight various applications of the proposed method, several types of well-studied tree models are selected to serve as representatives and then the corresponding solutions to geodesic distances on them are analytically derived. The result shows that our consequences are completely the same as previously reported ones, implying which our method is consistent.

(3) With the help of both our obtained results and effective resistance of electric networks, we also precisely determine the expressions of mean first-passage time ($MFPT$) for random walk on all trees considered here. The results mean that the underlying structures have a considerable influence on scale of both exact solutions of $MFPT$ and vertex number of tree in question.

This paper can be organized as follows. In Section 2, we will build up some helpful operations, including ($1,m$)-star-fractal operation, to generate our desired tree models, and introduce useful notations, for example, surjection and bijection, to develop our main results. And then, using a novel method proposed in Section 3, we derive exact solutions for geodesic distances on two classes of tree models constructed here. Additionally, we also obtain the closed-form solution for geodesic distance on a special case of the resulting tree models created by the ($1,m$)-star-fractal operation in a manner with respect to its own self-similarity. While the solution is the same as that derived using our novel method, such a manner can not be adequately employed under more general circumstance, suggesting that our method is more general. To show some other potential applications of our methods, in Section 4, we make use of connection between random walk and electrical network to analytically derive solutions for mean first-passage time on all the generated tree models. Related work is shown in Section 5. In conclusion, we close this paper in Section 6.

\begin{figure}
\centering
  \includegraphics[height=2.8cm]{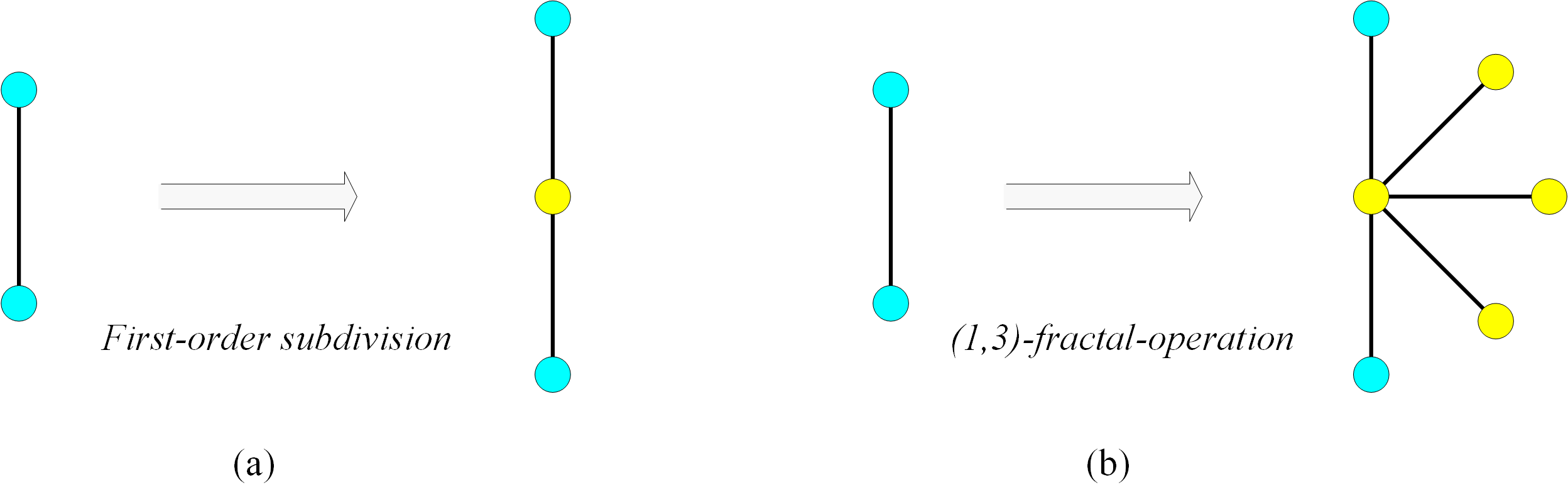}\\
{\small Fig.1. The diagrams of two types of operations on an edge. First-order subdivision on an edge is shown in panel (a) and the panel (b) shows ($1,m$)-fractal-operation on an edge where $m=3$.       }
\end{figure}

\section{Definitions and notations}

Here, we will recall some fundamental definitions and widely adopted notations from graph theory. It is conventional to let symbol $\mathcal{G}(\mathcal{V},\mathcal{E})$ denote as a graph where $\mathcal{V}$ and $\mathcal{E}$ are vertex set and edge set, respectively, and its corresponding order (vertex number) and size (edge number) are denoted by $|\mathcal{V}|$ and $|\mathcal{E}|$ where symbol $||$ represents the cardinality of a set. In the meantime, we denote by the notation $[1,n]$ an integer set which  consists precisely of those integers both no larger than $n$ and no less than $1$. For more details about other used not yet defined notions to see \cite{Bondy-2008}.

\textbf{Definition 1} Given an arbitrary graph $\mathcal{G}(\mathcal{V},\mathcal{E})$, if one inserts a new vertex $w$ to every edge $uv\in \mathcal{E}$ then the resulting graph, defined as $\mathcal{G'}_{1}(\mathcal{V'}_{1},\mathcal{E'}_{1})$, is called the \textbf{\emph{first-order subdivision graph }}of original graph $\mathcal{G}(\mathcal{V},\mathcal{E})$. To put this another way, the first-order subdivision graph can be equivalently obtained from graph $\mathcal{G}(\mathcal{V},\mathcal{E})$ by replacing every edge $uv\in \mathcal{E}$ by a unique path $uwv$ with length $2$ where internal vertex $w$ is in fact that inserted vertex. Henceforth, we regard such an operation on each edge of a graph as \textbf{\emph{first-order subdivision}}. It is worth noting that in this paper, we focus mainly on discussions about impact from first-order subdivision on geodesic distance of tree $\mathcal{T}(\mathcal{V},\mathcal{E})$ of significant interest. Here, Fig.1(a) illustrates the first-order subdivision on an edge.

For our purpose, it can immediately know by Def.1 that the first-order subdivision graph $\mathcal{G'}_{1}(\mathcal{V'}_{1},\mathcal{E'}_{1})$ holds on a couple of equations $|\mathcal{V'}_{1}|=|\mathcal{V}|+|\mathcal{E}|$ and $|\mathcal{E'}_{1}|=2|\mathcal{E}|$. After iteratively applying first-order subdivision on graph $\mathcal{G}(\mathcal{V},\mathcal{E})$ $t$ times, the order $|\mathcal{V'}_{t}|$ and size $|\mathcal{E'}_{t}|$ of the first-order subdivision graph $\mathcal{G'}_{t}(\mathcal{V'}_{t},\mathcal{E'}_{t})$ will follow a pair of equations

\begin{equation}\label{Section-2-0}
|\mathcal{V'}_{t}|=|\mathcal{V}|+(2^{t}-1)|\mathcal{E}|, \qquad |\mathcal{E'}_{t}|=2^{t}|\mathcal{E}|.
\end{equation}

\textbf{Definition 2} Given an arbitrary graph $\mathcal{G}(\mathcal{V},\mathcal{E})$, if one not only inserts a vertex $w$ to every edge $uv\in \mathcal{E}$ but also connects $m$ other new vertices $w_{i}$ ($i\in [1,m]$) to this newly inserted vertex $w$, then the resulting graph, referred to as $\mathcal{G^{\star}}_{1}(\mathcal{V^{\star}}_{1},\mathcal{E^{\star}}_{1})$, is called the\textbf{ \emph{(1,m)-star-fractal graph} }of original graph $\mathcal{G}(\mathcal{V},\mathcal{E})$. Equivalently speaking, such an operation can be manipulated on graph $\mathcal{G}(\mathcal{V},\mathcal{E})$ by directly inserting a star, where vertex $w$ is the center attached to $m$ leaves $w_{i}$, to every edge $uv\in \mathcal{E}$ and hence called\textbf{ \emph{(1,m)-star-fractal operation}}. It is obvious to say that the well known \textbf{\emph{T-fractal}} can be induced as a particular case of our ($1,m$)-star-fractal graph when parameter $m$ is supposed equal to $1$ \cite{Peng-2018}. As before, our aim is to study geodesic distance on each member of ($1,m$)-star-fractal tree family. An example as illustration of ($1,m$)-star-fractal operation on an edge is shown in Fig.1(b) where the newly generated star has $3$ leaves.

For brevity, with the help of Def.2, one can find out that ($1,m$)-star-fractal graph $\mathcal{G^{\star}}_{1}(\mathcal{V^{\star}}_{1},\mathcal{E^{\star}}_{1})$ has $|\mathcal{V^{\star}}_{1}|=|\mathcal{V}|+(1+m)|\mathcal{E}|$ vertices and $|\mathcal{E^{\star}}_{1}|=(2+m)|\mathcal{E}|$ edges. Similarly, for ($1,m$)-star-fractal graph $\mathcal{G^{\star}}_{t}(\mathcal{V^{\star}}_{t},\mathcal{E^{\star}}_{t})$, its vertex number $|\mathcal{V^{\star}}_{t}|$ and edge number $|\mathcal{E^{\star}}_{t}|$, respectively, obey

\begin{equation}\label{Section-2-1}
|\mathcal{V^{\star}}_{t}|=|\mathcal{V}|+((2+m)^{t}-1)|\mathcal{E}|, \qquad |\mathcal{E^{\star}}_{t}|=(2+m)^{t}|\mathcal{E}|.
\end{equation}

The statements above are mainly correlated with an arbitrary graph. However, below will introduce a helpful definition about vertex cover on a path in the jargon of graph theory that will be proved useful to derive our desired results in subsequent sections. This is in essence a kind of classification of vertices on a path. For more details to see \cite{Bondy-2008}, we here just discuss such a problem in the simplest and most fundamental situation.

\textbf{Definition 3} Given a path $\mathcal{P}$ with $n$ vertices, it is easy to see that these $n$ vertices can be divided into two disconnected vertex sets, without loss of generality, labeled as set $X=\{x_{1},x_{2},...,x_{\lceil n/2\rceil}\}$ and set $Y=\{y_{1},y_{2},...,y_{\lfloor n/2\rfloor}\}$. The above classification of vertices is in fact a bipartition of vertex set. More generally, vertices of set $X$ and vertices of set $Y$ can be alternatively arranged on the path $\mathcal{P}$ in an appropriate manner such that arbitrary vertex pair $x_{i}$ and $x_{j}$ is not connected directly by an edge and similarly for all pairs of vertices $y_{i}$ and $y_{j}$. This suggests that both vertex set $X$ and vertex set $Y$ are \textbf{\emph{vertex covers}} of path $\mathcal{P}$ where $|X|=|Y|+1$ when $n$ is odd and $|X|=|Y|$ otherwise. The both vertex sets will be alternatively employed to play a vital role in the process of building up our main results in the rest of this paper.

\textbf{Definition 4} Given two non-empty sets $X$ and $Y$, one can make a \textbf{\emph{mapping}} $f$ from $X$ to $Y$ such that for a provided element $x$ of set $X$ there must be a unique element $y$ in set $Y$ satisfying $f(x)=y$. This can be simply expressed in the following

$$\forall x\in X, \quad \exists y\in Y, \quad s.t.,\quad f:\; x\mapsto \; y$$
where the image set of set $X$ may be referred to as $f_{X\mapsto Y}=\{y|f(x)=y, \; y\in Y\}$. In this paper, we are interested in the below two kinds of mappings between sets $X$ and $Y$.

\emph{case 1} If both $|X|>|Y|$ and $|f_{X\mapsto Y}|=|Y|$ hold true, then this mapping $f$ is considered \textbf{\emph{surjection}}. Besides, for each element $y$ of image set $Y$, if there exist $n$ distinct pre-images $x_{i}$ ($i\in [1,n]$), i.e., $f^{-1}(y)=x_{i}$, then the surjection $f$ is considered \textbf{\emph{n-regular}}. It is clear to the eye that both sets $X$ and $Y$ follow $|X|=n |Y|$ when surjection $f$ is $n$-regular.

\emph{case 2} If the surjection $f$ under consideration holds $|X|=|Y|$ then it can be thought of as a \textbf{\emph{bijection}}, also called one-one mapping.

For convenience, the compound mapping between two mappings $f$ and $g$ can be expressed as $f\circ g$ mathematically.

So far, we have introduced some helpful definitions and notations used later. As stated above, the topic of this paper focuses principally on many discussions correlated to geodesic distance on tree models of significant interest. In addition, we will take useful advantage of a trivial property of tree that any pair of vertices $u$ and $v$ of a tree is connected by a unique path $\mathcal{P}_{uv}$. In fact, the length of such a path is indeed the geodesic distance $d_{uv}$ of this vertex pair. For simplicity and convenience, we abuse vertex pair $<u,v>$ to indicate a path whose endvertices are $u$ and $v$. As a result, the geodesic distance $\mathcal{S}$ of a graph $\mathcal{G}(\mathcal{V},\mathcal{E})$ can be written as

\begin{equation}\label{Section-2-2}
\mathcal{S}=\sum_{u,v\in\mathcal{V}}d_{uv}
\end{equation}
where vertex $u$ is distinct with $v$. An illustrative example for determining exact solution to geodesic distance on a tree is provided in Sec.1 of Supplementary Materials. It is natural to denote by $\langle\mathcal{S}\rangle=2\mathcal{S}/[|\mathcal{V}|(|\mathcal{V}|-1)]$ average geodesic distance on graph $\mathcal{G}(\mathcal{V},\mathcal{E})$.

\section{Main results}

We will in this section show main results, which are organized into several theorems, corollaries and applications in form, using the novel methods described later.

\textbf{Theorem 1} Given an arbitrary tree $\mathcal{T}(\mathcal{V},\mathcal{E})$, the exact solution for geodesic distance $\mathcal{S'}_{1}$ of its first-order subdivision tree $\mathcal{T'}_{1}(\mathcal{V'}_{1},\mathcal{E'}_{1})$ is

\begin{equation}\label{Section-3-1-0}
\mathcal{S'}_{1}=8\mathcal{S}-2|\mathcal{V}|(|\mathcal{V}|-1)
\end{equation}
in which $\mathcal{S}$ is a known expression to geodesic distance of tree $\mathcal{T}(\mathcal{V},\mathcal{E})$.

\textbf{\emph{Proof}} Consider an arbitrary tree $\mathcal{T}(\mathcal{V},\mathcal{E})$, first, assume that its geodesic distance is equal to $\mathcal{S}$. After applying first-order subdivision to each edge of tree $\mathcal{T}(\mathcal{V},\mathcal{E})$, the first-order subdivision tree $T'_{1}(\mathcal{V'}_{1},\mathcal{E'}_{1})$ will consist of two different types of vertices, without loss of generality, which are grouped into two disjoint vertex sets $X'$ and $Y'$. Set $X'$ contains the total old vertices of original tree $\mathcal{T}(\mathcal{V},\mathcal{E})$ and is in fact set $\mathcal{V}$, namely, $x'\in X'$ being the same as $x\in \mathcal{V}$. The other set $Y'$ is constituted by all the created vertices using the first-order subdivision. Obviously, $Y'=\mathcal{V'}_{1}-\mathcal{V}=\mathcal{V'}_{1}-X'$. In order to precisely calculate geodesic distance $\mathcal{S'}$, it is straightforward to compute three classes of geodesic distances, one for vertex pairs $<x'_{i},x'_{j}>$ of set $X'$, one for vertex pairs $<y'_{i},y'_{j}>$ of set $Y'$ as well as the latter for vertex pairs $<x'_{i},y'_{j}>$ between set $X'$ and set $Y'$. To this end, we will in turn accomplish these calculations according to complexity.

\emph{Case 1.1} For a given vertex pair $<x'_{i},x'_{j}>$ of set $X'$, there must be a bijection $f_{1}$ between set $X'$ and set $\mathcal{V}$ such that $f_{1}(<x'_{i},x'_{j}>)=<x_{i},x_{j}>$. Here vertices $x_{i}$ and $x_{j}$ are in set $\mathcal{V}$. In fact, such a bijection $f_{1}$ is self-mapping and so one can write

\begin{equation}\label{Section-3-1-1}
\mathcal{S'}_{1}(1)=2\mathcal{S}
\end{equation}
where $\mathcal{S'}_{1}(1)$ is the sum of distances of all possible vertex pairs in set $X'$. This is a consequence directly related to the intrinsic nature of first-order subdivision. The left tasks are to look for a reliable relation connecting the equations to be built to Eq.(\ref{Section-3-1-1}) just because of Eq.(\ref{Section-3-1-1})'s own simplicity.

\emph{Case 1.2} Compared to \emph{case 1.1}, there indeed exists a self-mapping between set $Y'$ and set $\mathcal{V'}-\mathcal{V}$. But, this has no help for addressing our problem of determining sum of geodesic distances on the second class of vertex pairs. Taking into account results in \emph{case 1.1 }known to us, the current goal is to build a connection to Eq.(\ref{Section-3-1-1}). More specifically, for a given pair of vertices $<y'_{i},y'_{j}>$ belonging to set $Y'$ as shown in Fig.2, we may find out a bijection $f_{2}$ between set $Y'$ and set $X'$ that satisfies our requirement by means of statements in both Def.3 and Def.4. To see why this is so, let us pay attention on the both disjoint vertex sets $X'$ and $Y'$. For arbitrary vertex pair $<y'_{i},y'_{j}>$ of set $Y'$, there must be a unique path $\mathcal{P}_{y'_{i}y'_{j}}$ connecting vertices $y'_{i}$ and $y'_{j}$ which is lighted in red as shown in Fig.2 (online). At the same time, one certainly derives an extension $\mathcal{P}_{x'_{i}x'_{j}}$ from path $\mathcal{P}_{y'_{i}y'_{j}}$ by jointing two edges $x'_{i}y'_{i}$ and $x'_{j}y'_{j}$ (lighted in blue as shown in Fig.2), indicating which there is a mapping $f^{*}$ between vertex pairs $<y'_{i},y'_{j}>$ and $<x'_{i},x'_{j}>$. Meanwhile, it is not hard to turn out mapping $f^{*}$ to be bijection according to intrinsic properties among bijection $f_{1}$, first-order subdivision and tree itself. Thus, the compound mapping between the candidate $f^{*}$ and bijection $f_{1}$ may be designed as our desired bijection $f_{2}$, i.e., $f_{2}=f_{1}\circ f^{*}$. We then have

\begin{equation}\label{Section-3-1-2}
\mathcal{S'}_{1}(2)=\mathcal{S'}_{1}(1)-2\frac{|\mathcal{V}|(|\mathcal{V}|-1)}{2}
\end{equation}
in which $\mathcal{S'}_{1}(2)$ is geodesic distance of all possible vertex pairs in set $Y'$.

\begin{figure}
\centering
  \includegraphics[height=5cm]{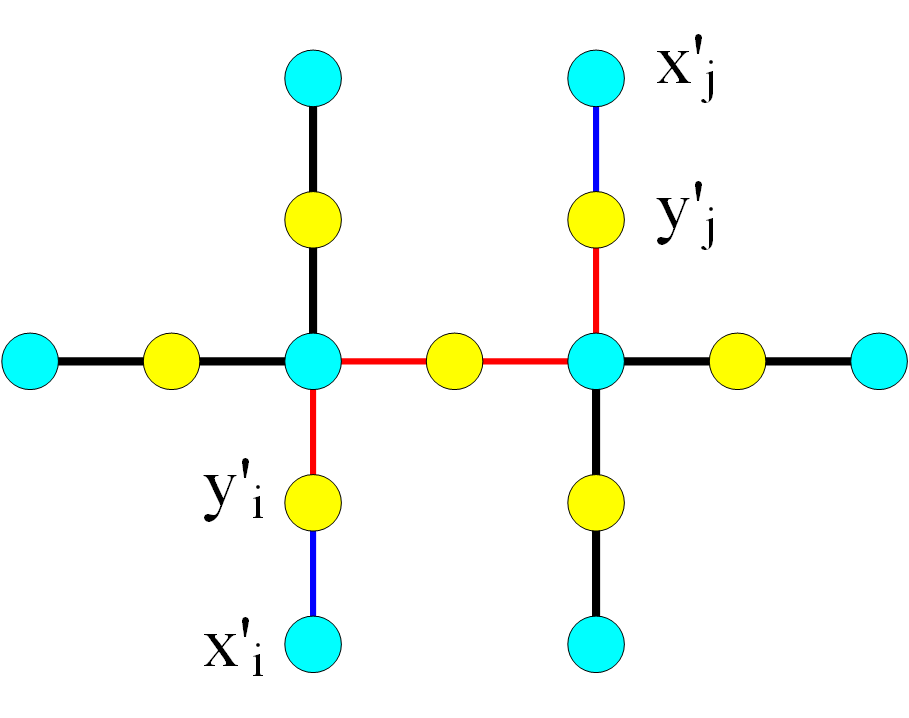}\\
{\small Fig.2. The illustration of proof in \emph{Case 1.2}.       }
\end{figure}

\emph{Case 1.3} The remainder of our problem to answer is to capture the expression of geodesic distances $\mathcal{S'}_{1}(3)$ for all possible vertex pairs $<x'_{i},y'_{j}>$ whose endvertices are from different sets, i.e., set $X'$ and set $Y'$. Along the research line of \emph{case 1.2}, for vertex pair $x'_{i}$ and $y'_{j}$, there is also a unique path $\mathcal{P}_{x'_{i}y'_{j}}$ which can be reduced to another path $\mathcal{P}_{y'_{i}y'_{j}}$ by deleting an additional edge $x'_{i}y'_{i}$ based on a similar mapping $f^{**}$ to mapping $f^{*}$ in \emph{case 1.2}. Therefore, we may take a mapping $f_{3}=f_{1}\circ f^{**}$ that must be a surjection between vertex pairs $<x'_{i},y'_{j}>$ and $<x'_{i},x'_{j}>$ but not bijection. One of reasons for this is that another path $\mathcal{P}_{y'_{i}x'_{j}}$ may be induced as the path $\mathcal{P}_{y'_{i}y'_{j}}$ by removing edge $y'_{j}x'_{j}$ as well. Armed with the both cases, there are two distinct pre-images under surjection $f_{3}$, i.e., $<x'_{i},y'_{j}>$ and $<y'_{i},x'_{j}>$, such that $f^{-1}_{3}(<y'_{i},y'_{j}>)=<x'_{i},y'_{j}>=<y'_{i},x'_{j}>$. By Def.4, such a surjection is in principle 2-regular and hence the exact solution for geodesic distance $\mathcal{S'}_{1}(3)$ obeys

\begin{equation}\label{Section-3-1-3}
\mathcal{S'}_{1}(3)=2\mathcal{S'}_{1}(2)+|\mathcal{\mathcal{V}}|(|\mathcal{V}|-1).
\end{equation}

With the declarations in \emph{cases 1.1-1.3}, Eqs.(\ref{Section-3-1-1})-(\ref{Section-3-1-3}) together produce an exact solution for $\mathcal{S'}_{1}$ completely equivalent to that of Eq.(\ref{Section-3-1-0}). This completes our proof. The detailed calculation of geodesic distance $\mathcal{S'}_{1}$ on tree shown in Fig.2 is provided in Sec.2 of Supplementary Materials for purpose of illustration.

On the basis of Eq.(\ref{Section-3-1-0}), we directly give the solution for geodesic distance $\mathcal{S'}_{t}$ on the first-order subdivision tree $\mathcal{T'}_{t}(\mathcal{V'}_{t},\mathcal{E'}_{t})$ with omitting detailed computations, as follows

\emph{Corollary 1} After $t$ time steps, the solution of geodesic distance $\mathcal{S'}_{t}$ on the first-order subdivision tree $\mathcal{T'}_{t}(\mathcal{V'}_{t},\mathcal{E'}_{t})$ will follow

\begin{equation}\label{Section-3-c-10}
\mathcal{S'}_{t}=8^{t}\mathcal{S}-\frac{1}{3}(2^{3t}-2^{t})(|\mathcal{V}|-1)+(2^{2t-1}-2^{3t-1})(|\mathcal{V}|-1)^{2}.
\end{equation}

\textbf{Application 1} After $t$ time steps, the solution of average geodesic distance $\langle\mathcal{S'}_{t}\rangle$ on the first-order subdivision tree $\mathcal{T'}_{t}(\mathcal{V'}_{t},\mathcal{E'}_{t})$ will follow

\begin{equation}\label{Section-3-a-10}
\begin{aligned}\langle\mathcal{S'}_{t}\rangle&=\frac{\mathcal{S'}_{t}}{|\mathcal{V'}_{t}|(|\mathcal{V'}_{t}|-1)/2}\\
&\approx \frac{2^{t+1}\mathcal{S}}{(|\mathcal{V}|-1)^{2}}-\frac{2^{t+1}}{3(|\mathcal{V}|-1)}+1-2^{t}\\
&=O(|\mathcal{V'}_{t}|)
\end{aligned}
\end{equation}
where we have made use of result in Eq.(\ref{Section-2-0}).

As known, the most special member of tree family is an edge. If let the seed of the first-order subdivision tree $\mathcal{T'}_{t}(\mathcal{V'}_{t},\mathcal{E'}_{t})$ be an edge connecting a couple of vertices, then model $\mathcal{T'}_{t}(\mathcal{V'}_{t},\mathcal{E'}_{t})$ will become a path $\mathcal{P'}_{t}(\mathcal{V'},\mathcal{E'})$ with $2^{t}+1$ vertices. Therefore, we are able to state the following corollary according to Eq.(\ref{Section-3-c-10}).

\emph{Corollary 2}  After $t$ time steps, the solution of geodesic distance $\mathcal{S}(t)$ on the first-order subdivision path $\mathcal{P'}_{t}(\mathcal{V'},\mathcal{E'})$ will follow

\begin{equation}\label{Section-3-c-20}
\mathcal{S}(t)=\frac{(2^{t-1}+1)2^{t}(2^{t}+1)}{3}.
\end{equation}

Equivalently, it is not hard to capture the solution of geodesic distance $\mathcal{S}(t)$ on path of such type in the most general manner, i.e. enumeration method, as follows

\begin{equation}\label{Section-3-c-21}
\begin{aligned}\mathcal{S}(t)&=\sum_{i=1}^{2^{t}}\sum_{j=1}^{2^{t}+1-i}j\\
&=\sum_{i=1}^{2^{t}}\frac{(2^{t}+1-i)(2^{t}+2-i)}{2}.
\end{aligned}
\end{equation}

Clearly, Eqs.(\ref{Section-3-c-20}) and (\ref{Section-3-c-21}) here give a concise proof for one combinatorial identity $$\sum_{i=1}^{2^{t}}\sum_{j=1}^{2^{t}+1-i}j=\frac{(2^{t-1}+1)2^{t}(2^{t}+1)}{3}$$
based on computation of geodesic distance $\mathcal{S}(t)$ on path with length $2^{t}$. Put this further, one can easily derive the average geodesic distance $\langle\mathcal{S}(t)\rangle$ on the first-order subdivision path $\mathcal{P'}_{t}(\mathcal{V'}_{t},\mathcal{E'}_{t})$ from Eq.(\ref{Section-3-c-20}).

\textbf{Application 2} After $t$ time steps, the solution of average geodesic distance $\langle\mathcal{S}(t)\rangle$ on the first-order subdivision path $\mathcal{P'}_{t}(\mathcal{V'},\mathcal{E'})$ will follow

\begin{equation}\label{Section-3-a-20}
\langle\mathcal{S}(t)\rangle\approx \frac{2^{t}+3}{3}.
\end{equation}

Taking into account results between Eq.(\ref{Section-3-a-10}) and Eq.(\ref{Section-3-a-20}), we can immediately capture the below theorem which says an interesting phenomenon about influence from the first-order subdivision on average geodesic distance on the first-order subdivision tree $\mathcal{T'}_{t}(\mathcal{V'}_{t},\mathcal{E'}_{t})$.

\textbf{Theorem 2}  Given an arbitrary tree $\mathcal{T}(\mathcal{V},\mathcal{E})$, the solution for average geodesic distance $\langle\mathcal{S'}_{t}\rangle$ on its first-order subdivision tree $\mathcal{T'}_{t}(\mathcal{V'}_{t},\mathcal{E'}_{t})$ will follow

\begin{equation}\label{Section-3-t-20}
\langle\mathcal{S'}_{t}\rangle \approx|\mathcal{V'}_{t}|^{\gamma'}= O(D'_{t})
\end{equation}
where exponent $\gamma'=1$ in the large graph size limit and we denote by $D'_{t}$ diameter of tree $\mathcal{T'}_{t}(\mathcal{V'}_{t},\mathcal{E'}_{t})$ \footnote[1]{In the language of graph theory, the diameter $D$ is the maximum over geodesic distances of all vertex pairs in a graph $\mathcal{G}(\mathcal{V},\mathcal{E})$. }.

By far, the first-order subdivision trees $\mathcal{T'}_{t}(\mathcal{V'}_{t},\mathcal{E'}_{t})$ which are generated by an arbitrary tree considered as a seed all share some features in common including: (i) Each vertex added by first-order subdivision at any time step $t_{i}$ ($1\leq t_{i}\leq t$) into models $\mathcal{T'}_{t}(\mathcal{V'}_{t},\mathcal{E'}_{t})$ has degree $2$, (ii) The total number of leaves of models $\mathcal{T'}_{t}(\mathcal{V'}_{t},\mathcal{E'}_{t})$ keeps unchanged \cite{Ma-2019}, and (iii) All models $\mathcal{T'}_{t}(\mathcal{V'}_{t},\mathcal{E'}_{t})$ exhibit homogeneous topological structure. Besides that, hereafter, the rest of this section will discuss another type of tree models $\mathcal{T}^{\star}_{t}(\mathcal{V}^{\star}_{t},\mathcal{E}^{\star}_{t})$ with inheterogeneous topological structure. These models not only have some similar properties to models $\mathcal{T'}_{t}(\mathcal{V'}_{t},\mathcal{E'}_{t})$ but also inherit some intriguing characters from ($1,m$)-star-fractal operation unseen in models $\mathcal{T'}_{t}(\mathcal{V'}_{t},\mathcal{E'}_{t})$, particularly, fractal property.

At first, let us start from investigating the simplest form of tree models $\mathcal{T}^{\star}_{t}(\mathcal{V}^{\star}_{t},\mathcal{E}^{\star}_{t})$ by introducing theorem 3.

\textbf{Theorem 3} Given an arbitrary tree $\mathcal{T}(\mathcal{V},\mathcal{E})$, the exact solution for geodesic distance $\mathcal{S}^{\star}_{1}$ of its ($1,m$)-star-fractal tree $\mathcal{T}^{\star}_{1}(\mathcal{V}^{\star}_{1},\mathcal{E}^{\star}_{1})$ is

\begin{equation}\label{Section-3-t-30}
\mathcal{S^{\star}}_{1}=2(m+2)^{2}\mathcal{S}-(m+2)(|\mathcal{V}|-1)(m+|\mathcal{V}|)
\end{equation}
in which $\mathcal{S}$ is a known expression to geodesic distance on tree $\mathcal{T}(\mathcal{V},\mathcal{E})$.

\emph{\textbf{Proof}} By Def.2, there also are two different groups of vertices, i.e., vertex sets $X'=\mathcal{V}$ and $Y'=\mathcal{V^{\star}}-\mathcal{V}$, and hence we still need to consider three kinds of contributions to computation of geodesic distance $\mathcal{S^{\star}}_{1}$ as discussed in the development of theorem 1. With the light shed by developing theorem 1, we will make use of more fine-grained classification method than that employed in the proof of theorem 1 to resolve this computations. Considering that the first-order subdivision is a special case of the ($1,m$)-star-fractal operation, we can directly use some existing results from the proof of theorem 1 to consolidate theorem 3 without specific descriptions. Meantime, some previously adopted notations keep unchanged. Besides, the vertex set $Y'$ should be partitioned into subsets $Y'_{1}$ and $Y'_{2}$ such that $Y'_{1}$ is constituted by the total leaf vertices of each newly inserted star and $Y'_{2}$ contains the central vertex of each newly created star. To distinguish newborn vertices between sets $Y'_{1}$ and $Y'_{2}$, each vertex in set $Y'_{2}$ remains marked $y'_{i}$ as above and then we label each vertex in set $Y'_{1}$ by $y''_{i}$. From now on, let us begin with clarifying the correctness of Eq.(\ref{Section-3-t-30}).

\emph{Case 2.1} For an arbitrary vertex pair $<x'_{i},x'_{j}>$ of set $X'$, under bijection $f_{1}$, the geodesic distance $\mathcal{S^{\star}}_{1}(1)$ on such type of vertex pairs complies to

\begin{equation}\label{Section-3-t-31}
\mathcal{S^{\star}}_{1}(1)=2\mathcal{S}.
\end{equation}

\emph{Case 2.2} For an arbitrary vertex pair $<y'_{i},y'_{j}>$ of set $Y'_{2}$, this surjection $f_{2}=f_{1}\circ f^{*}$ will make  geodesic distance $\mathcal{S^{\star}}_{1}(2)$ on such kind of vertex pairs satisfy

\begin{equation}\label{Section-3-t-32}
\mathcal{S^{\star}}_{1}(2)=\mathcal{S^{\star}}_{1}(1)-2\frac{|\mathcal{V}|(|\mathcal{V}|-1)}{2}.
\end{equation}

\emph{Case 2.3} For an arbitrary vertex pair $<x'_{i},y'_{j}>$ or $<y'_{i},x'_{j}>$, this surjection $f_{3}=f_{1}\circ f^{**}$ will guarantee geodesic distance $\mathcal{S^{\star}}_{1}(3)$ on such class of vertex pairs to obey

\begin{equation}\label{Section-3-t-33}
\mathcal{S^{\star}}_{1}(3)=2\mathcal{S^{\star}}_{1}(2)+|\mathcal{V}|(|\mathcal{V}|-1).
\end{equation}

\emph{Case 2.4} There must be $|\mathcal{V}|-1$ new stars introduced into original tree $\mathcal{T}(\mathcal{V},\mathcal{E})$ by means of ($1,m$)-star-fractal operation. Here we just capture geodesic distances on an arbitrary pair of leaf vertices within the same star, i.e., vertex pair $<y''_{oi},y''_{oj}>$ where the first subscript $o$ represents the central vertex of star attached to the both leaf vertices. Therefore the geodesic distance $\mathcal{S^{\star}}_{1}(4)$ on all possible leaf vertex pairs of this type follows

\begin{equation}\label{Section-3-t-34}
\mathcal{S^{\star}}_{1}(4)=(|\mathcal{V}|-1)\left(2\frac{m(m-1)}{2}\right).
\end{equation}

\emph{Case 2.5} We here discuss geodesic distance $\mathcal{S^{\star}}_{1}(5)$ on all possible leaf vertex pairs $<y''_{ui},y''_{vj}>$ in which two vertices come from different stars. As before, in order to accomplish this task, we have to choose a fresh mapping $f^{4*}$ that bridges vertex pairs $<y''_{ui},y''_{vj}>$ to vertex pair $<y'_{u},y'_{v}>$. Here both vertices $y'_{u}$ and $y'_{v}$ are, respectively, the central of stars to which vertices $y''_{ui}$ and $y''_{vj}$ belong. This anticipated mapping $f^{4*}$ will in essence connect any two stars and can be timely verified to be an $m^{2}$-regular surjection. And then, using bijection $f_{2}$ introduced in \emph{case 2.2}, we can generate a satisfactory surjection $f_{4}=f_{2}\circ f^{4*}$ in time and obtain

\begin{equation}\label{Section-3-t-35}
\mathcal{S^{\star}}_{1}(5)=m^{2}\mathcal{S^{\star}}_{1}(2)+\sum_{i=1}^{|\mathcal{V}|-2}2m^{2}(|\mathcal{V}|-1-i).
\end{equation}

\emph{Case 2.6} At the moment, let us pay attention to computation of geodesic distance $\mathcal{S^{\star}}_{1}(6)$ on all possible vertex pairs $<x'_{i},y''_{oj}>$ where $x'_{i}\in X'$ and $y''_{oj}\in Y'_{1}$. Note that, in some cases, two subscripts $i$ and $j$ can be equal. Considering such a vertex pair $<x'_{i},y''_{oj}>$ carefully, the first task is to find out a surjection $f^{5*}$ such that $f^{5*}(<x'_{i},y''_{oj}>)=<x'_{i},y'_{o}>$ where vertex $y'_{o}$ is the central vertex of star including vertex $y'_{oj}$. And then, using well proposed surjection $f_{3}$ in\emph{ case 2.3}, we can create an acceptable surjection $f_{5}=f_{3}\circ f^{5*}$ connecting vertex pair $<x'_{i},y''_{oj}>$ with $<x'_{i},y'_{o}>$ and so the solution for geodesic distance $\mathcal{S^{\star}}_{1}(6)$ is

\begin{equation}\label{Section-3-t-36}
\mathcal{S^{\star}}_{1}(6)=m\mathcal{S^{\star}}_{1}(3)+m|\mathcal{V}|(|\mathcal{V}|-1).
\end{equation}

\emph{Case 2.7} By now, we have successfully achieved the entire computations of geodesic distances on vertex pairs in which one vertex is from set $X'$ and the other belongs to set $Y'$. The issue to answer is to count geodesic distance $\mathcal{S^{\star}}_{1}(7)$ on vertex pairs where one vertex is selected from set $Y'_{1}$ and the other from set $Y'_{2}$. With the terminologies mentioned above, such a vertex pair can be thought of as $<y'_{i},y''_{oj}>$ in which it is possible that subscript $i$ is the same as $o$ when leaf vertex $y''_{oj}$ and central vertex $y'_{i}$ are in a common star. For all vertex pairs $<y'_{i},y''_{oj}>$ with $i\neq o$, it is natural to construct a $2m$-regular surjection $f^{6*}$ projecting vertex pair $<y'_{i},y''_{oj}>$ to $<y'_{i},y'_{o}>$, and then blurring new surjection $f^{6*}$ with bijection $f_{2}$ in \emph{case 2.2} provides us with an expectant surjection $f_{6}=f_{2}\circ f^{6*}$ that is able to be what we want. On the other hand, when $i$ is equal to $o$, the vertex pair $<y'_{i},y''_{ij}>$ will be mapped onto an identical vertex $y'_{i}$ under $m$-regular surjection $f^{7*}$. Through the descriptions here, a concise expression of geodesic distance $\mathcal{S^{\star}}_{1}(7)$ can be expressed as

\begin{equation}\label{Section-3-t-37}
\mathcal{S^{\star}}_{1}(7)=2m\mathcal{S^{\star}}_{1}(2)+m(|\mathcal{V}|-1)^{2}.
\end{equation}

Plugging Eqs.(\ref{Section-3-t-31})-(\ref{Section-3-t-37}) into this summarized expression $\mathcal{S^{\star}}_{1}=\sum_{i=1}^{7}\mathcal{S^{\star}}_{1}(i)$ and implementing some basic arithmetics together output the desirable result as said in Eq.(\ref{Section-3-t-30}). This suggests that theorem 3 is sound. As previously, we provide an example in Sec.3 of Supplementary Materials in order to show some details involved in developing theorem 3.

Similarly, we can immediately capture the solutions for geodesic distance $\mathcal{S^{\star}}_{t}$ and average geodesic distance $\langle\mathcal{S^{\star}}_{t}\rangle$ on the ($1,m$)-star-fractal tree $\mathcal{T^{\star}}_{t}(\mathcal{V^{\star}}_{t},\mathcal{E^{\star}}_{t})$, separately, which are stated in the next corollary and application.

\emph{Corollary 3} After $t$ time steps, the solution for geodesic distance $\mathcal{S^{\star}}_{t}$ on the ($1,m$)-star-fractal tree $\mathcal{T^{\star}}_{t}(\mathcal{V^{\star}}_{t},\mathcal{E^{\star}}_{t})$ will follow

\begin{equation}\label{Section-3-c-30}
\begin{aligned}\mathcal{S^{\star}}_{t}&=2^{t}(m+2)^{2t}\mathcal{S}+(2^{t}-1)(m+2)^{2t-1}(|\mathcal{V}|^{2}-2|\mathcal{V}|-1)
\\
&-\frac{(m+1)(|\mathcal{V}|-1)}{2}\times\frac{2^{t+1}(m+2)^{2t}-2(m+2)^{t}}{2(m+2)-1}
\end{aligned}
\end{equation}
in which $\mathcal{S}$ is a known expression of geodesic distance on tree $\mathcal{T}(\mathcal{V},\mathcal{E})$.

\textbf{Application 3} After $t$ time steps, the solution for average geodesic distance $\langle\mathcal{S^{\star}}_{t}\rangle$ on the ($1,m$)-star-fractal tree $\mathcal{T^{\star}}_{t}(\mathcal{V^{\star}}_{t},\mathcal{E^{\star}}_{t})$ will follow

\begin{equation}\label{Section-3-a-30}
\langle\mathcal{S^{\star}}_{t}\rangle\approx\frac{2^{t+1}\mathcal{S}}{(|\mathcal{V}|-1)^{2}}-\frac{(m+1)2^{t+1}}{(2m+3)(|\mathcal{V}|-1)}+1-2^{t}
\end{equation}
which is completely consistent with simulation results as plotted in Fig.3.

\begin{figure}
\centering
  \includegraphics[height=6cm]{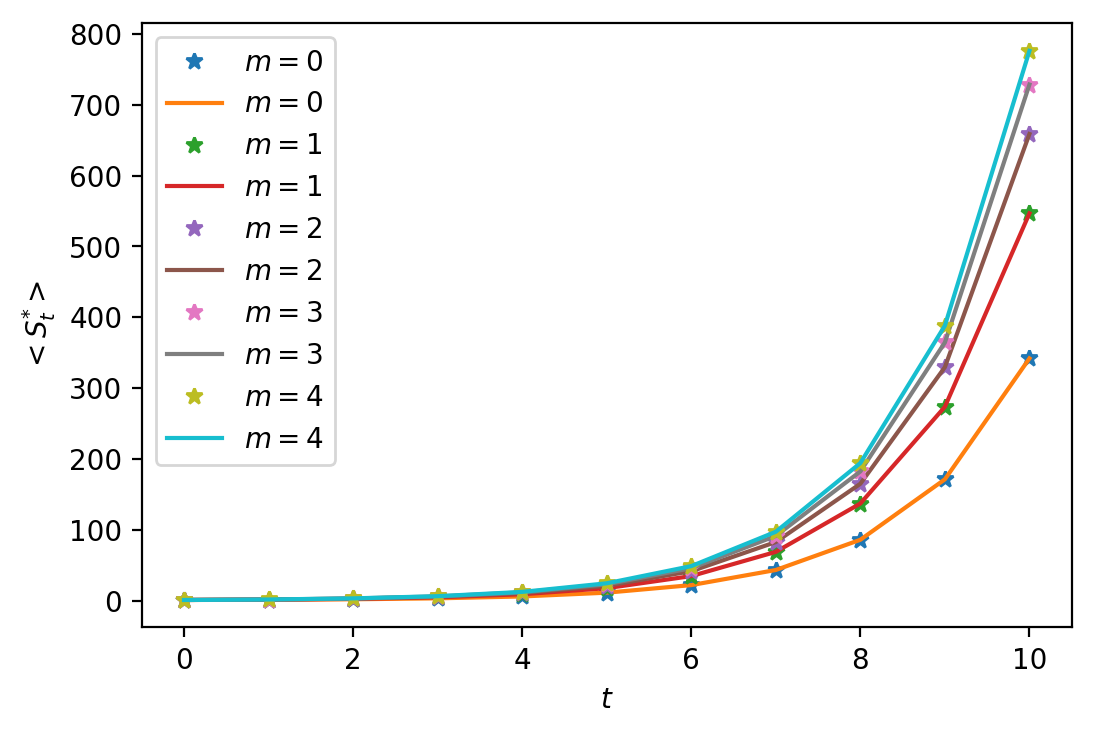}\\
{\small Fig.3. The diagram of average geodesic distance $\langle\mathcal{S^{\star}}_{t}\rangle$ on the ($1,m$)-star-fractal tree $\mathcal{T^{\star}}_{t}(\mathcal{V^{\star}}_{t},\mathcal{E^{\star}}_{t})$ where parameter $m$ is supposed equal to $0,1,2,3,4$, separately.       }
\end{figure}

As described before, there are in fact significant differences between the first-order subdivision and the ($1,m$)-star-fractal operation attributed to topological structures of these models generated by the both operations. Some of them will be reported in the rest of this paper in detail. Nevertheless, it can be quite evident that the two distinct kinds of tree models built by the two operations have similar expression of average geodesic distance, see Eq.(\ref{Section-3-a-10}) and Eq.(\ref{Section-3-a-30}). In another word, Eq.(\ref{Section-3-a-10}) can be regarded as a special case of results told by Eq.(\ref{Section-3-a-30}) when parameter $m$ is supposed equal to zero. This indirectly implies that the two operations above share similar function on some topological structural indices of generated tree models, at least on average geodesic distance.

\begin{figure*}
\centering
  \includegraphics[height=6.5cm]{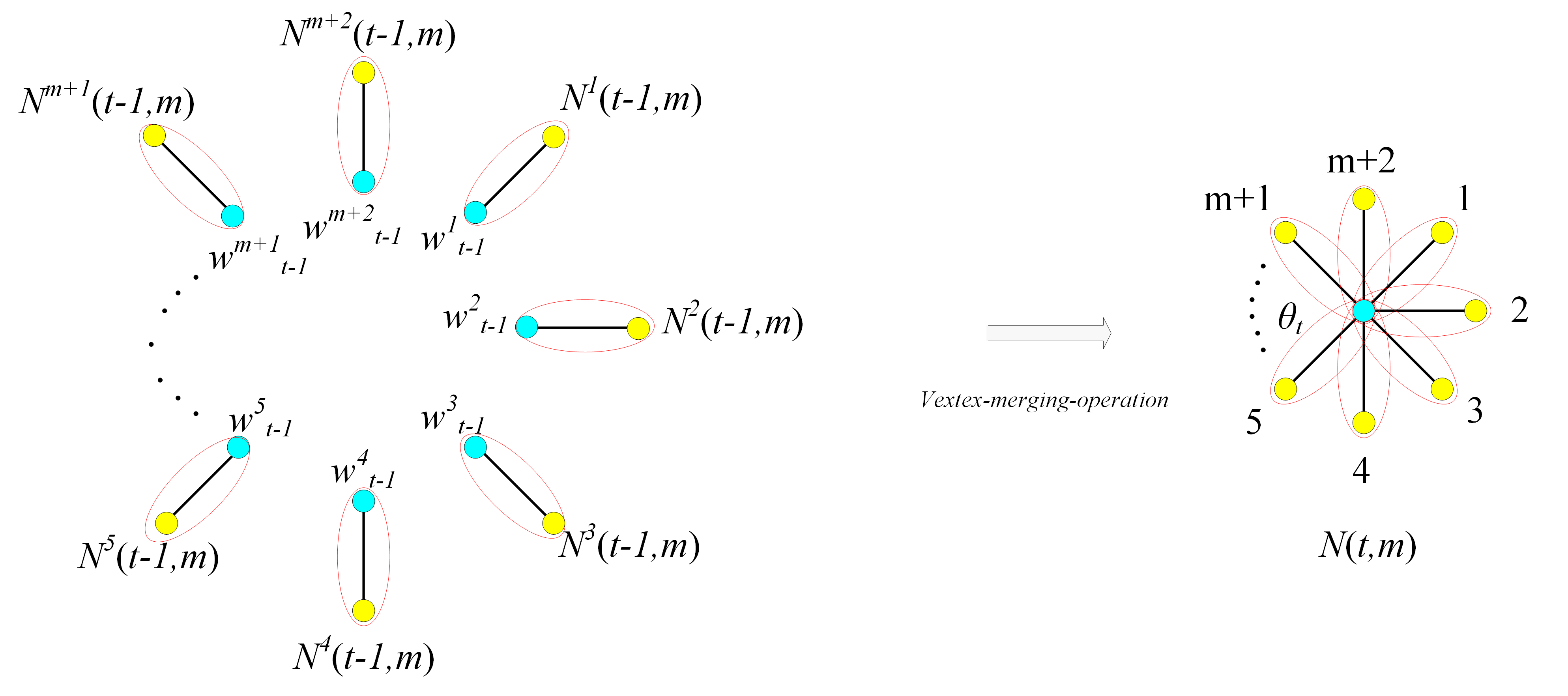}\\
{\small Fig.4. The diagram of tree model $N(t,m)$ which is in essence constructed by $m+2$ models $N(t-1,m)$ using vertex-merging-operation. Each small component of the left panel is an abstract representative of model $N(t-1,m)$. To be more concrete, we refer reader to the next figure, which describes T-graph($t$), for more details about reconstruction of such kind.}
\end{figure*}

Till now, the two families of tree models are established based on an arbitrary tree $\mathcal{T}(\mathcal{V},\mathcal{E})$ using the first-order subdivision and the ($1,m$)-star-fractal operation, respectively. More generally, our results answer how to analytically determine an exact solution for geodesic distance on tree models of such kinds. Therefore, the consequence published in \cite{Yuan-2010} can be viewed as a special example of our results in which a single edge is selected as the seed of growth tree models. To keep our work self-contained, we still study geodesic distance on such type of tree models \cite{Yuan-2010} and the corresponding precise expression is shown in corollary 4.

\emph{Corollary 4} After $t$ time steps, the solution of geodesic distance $\mathcal{S}(t,m)$ on tree model $N(t,m)$ will follow

\begin{equation}\label{Section-3-c-40}
\begin{aligned}\mathcal{S}(t,m)&=\frac{(m^{2}+2m+1)2^{t}+2m+3}{2m+3}(m+2)^{2t-1}\\
&+(m+2)^{t}-\frac{(m+2)^{2}}{2m+3}(m+2)^{t-1}
\end{aligned}.
\end{equation}
Here, tree model $N(t,m)$ is a special case of ($1,m$)-star-fractal tree $\mathcal{T^{\star}}_{t}(\mathcal{V^{\star}}_{t},\mathcal{E^{\star}}_{t})$ in which the selected seed is no longer an arbitrary tree but a single edge.

Tree models $N(t,m)$ of such type in fact have been in-depth studied in many published papers because they show some interesting structural features, such as, factual phenomena \cite{Yuan-2010}. Although Eq.(\ref{Section-3-c-40}) can be easily proved by plugging the initial conditions $\mathcal{S}_{0}=1$ and $|\mathcal{V}_{0}|=2$ into Eq.(\ref{Section-3-c-30}), we will turn out Eq.(\ref{Section-3-c-40}) to be correct in another fashion based on self-similarity displayed by tree models $N(t,m)$. The reasons why we employ method on the basis of self-similar topological structure have twofold. The one is that self-similarity is one of most prevailing topological structures of networked models in nature and real-life world \cite{Chao-2005}-\cite{Fei-2018}. The other is to highlight convenience of our novel methods addressed above in comparison with the commonly-used method which we will show below.

\textbf{\emph{Proof }} To smoothly develop the proof for Eq.(\ref{Section-3-c-40}), we have to describe the development process of tree models $N(t,m)$ by utilizing another reconstruction method, shown in Fig.4, where we denote by $\theta_{t}$ the center (indigo online) that can be obtained by vertex-merging-operation among the external vertices $\omega^{i}_{t-1}$\footnote[2]{In practice, the external vertex may be either of both vertices connected by the original edge as seed.} of branches $N^{i}(t-1,m)$ ($i\in [1,m+2]$). Indeed, this reconstruction method shows self-similar structure of tree models $N(t,m)$ and further allows us to calculate the exact solution of geodesic distance $\mathcal{S}(t,m)$ analytically, as below

\begin{equation}\label{Section-3-c-41}
\mathcal{S}(t,m)=(m+2)\mathcal{S}(t-1,m)+\Omega_{t,m}
\end{equation}
where symbol $\Omega_{t,m}$ represents the total sum of geodesic distance of an arbitrary pair of vertices from two different branches $N^{i}(t-1,m)$ ($i\in [1,m+2]$). Obviously, Eq.(\ref{Section-3-c-41}) can be reorganized in an iterative calculation way as follows

\begin{equation}\label{Section-3-c-42}
\mathcal{S}(t,m)=(m+2)^{t}\mathcal{S}(0,m)+\sum_{i=0}^{t-1}(m+2)^{i}\Omega_{t-i,m}
\end{equation}
where $\mathcal{S}(0,m)$ is the geodesic distance of two vertices connected by the original edge as a seed and in fact equals $1$.

To capture the closed-form solution of Eq.(\ref{Section-3-c-42}), the left issue is to answer the expressions of $\Omega_{t-i,m}$ ($i\in [0,t-1]$). Now, we, by definition, write

 \begin{equation}\label{Section-3-c-43}
\Omega_{t,m}=\sum_{1\leq i<j\leq m+2}\Omega^{ij}_{t,m}=\frac{(m+1)(m+2)}{2}\Omega^{12}_{t,m}
\end{equation}
in which we have made use of self-similar structure among branches $N^{i}(t-1,m)$ ($i\in [1,m+2]$). As before, we can by definition obtain

 \begin{equation}\label{Section-3-c-44}
\begin{aligned}\Omega^{12}_{t,m}&=\sum_{\begin{aligned}&v\in \mathcal{V}^{1}_{t-1},v\neq \omega^{1}_{t-1}( or \neq \theta_{t})\\
&u\in \mathcal{V}^{2}_{t-1},u\neq \omega^{2}_{t-1}( or \neq \theta_{t})
\end{aligned}}d^{vu}_{t}\\
&=\sum_{\begin{aligned}&v\in \mathcal{V}^{1}_{t-1},v\neq \omega^{1}_{t-1}( or \neq \theta_{t})\\
&u\in \mathcal{V}^{2}_{t-1},u\neq \omega^{2}_{t-1}( or \neq \theta_{t})
\end{aligned}}(d^{v\theta_{t}}_{t}+d^{\theta_{t}u}_{t})\\
&=2(|\mathcal{V}_{t-1}|-1)\Theta_{t-1}
\end{aligned}
\end{equation}
here we define $\Theta_{t-1}$ as the total sum of geodesic distances between the external vertex $ \omega^{1}_{t-1}$ and vertices $v\neq \omega^{1}_{t-1}( or \neq \theta_{t})$ of branch $N^{1}(t-1,m)$. We take useful advantage of self-similar structure between branches $N^{1}(t-1,m)$ and $N^{2}(t-1,m)$ again. Analogously, $\Theta_{t-1}$ can be written as

 \begin{equation}\label{Section-3-c-45}
\begin{aligned}\Theta_{t-1}&=\sum_{v\in \mathcal{V}_{t-1},v\neq \omega_{t-1}}d^{v\omega_{t-1}}_{t-1}\\
&=\Theta_{t-2}+\sum_{j\in[2,m+2]}\sum_{i\in \mathcal{V}^{j}_{t-2},i\neq \omega^{j}_{t-2}}(d^{i\omega^{j}_{t-2}}+D_{t-2})\\
&=(m+2)\Theta_{t-2}+(m+1)(|\mathcal{V}_{t-2}|-1)D_{t-2}
\end{aligned}
\end{equation}
where symbol $D_{t-2}$ is the diameter of tree model $N(t-2,m)$ and self-similarity among branches $N^{i}(t-2,m)$ ($i\in [2,m+2]$) is again used for simplicity. With the similar calculation to Eq.(\ref{Section-3-c-42}), the closed-form solution of $\Theta_{t-1}$ can follow

 \begin{equation}\label{Section-3-c-46}
\Theta_{t-1}=(m+1)\sum_{i=0}^{t-2}(m+2)^{i}(|\mathcal{V}_{t-2-i}|-1)D_{t-2-i}+(m+2)^{t-1}\Theta_{0}.
\end{equation}

Substituting both initial conditions $\Theta_{0}=1$ and $D_{t}=2^{t}$ into Eq.(\ref{Section-3-c-46}) yields

 \begin{equation}\label{Section-3-c-47}
\Theta_{t-1}=(m+2)^{t-1}+(m+1)(m+2)^{t-2}(2^{t-1}-1).
\end{equation}

Armed with Eqs.(\ref{Section-3-c-42})-(\ref{Section-3-c-47}), the exact solution of geodesic distance on tree models $N(t,m)$ may obey

 \begin{equation}\label{Section-3-c-48}
\begin{aligned}\mathcal{S}(t,m)&=(m+2)^{t}+(m+2)^{2t-1}-(m+2)^{t-1}\\
&+\frac{(m+1)^{2}}{2m+3}[2^{t}(m+2)^{2t-1}-(m+2)^{t-1}]
\end{aligned}.
\end{equation}

By some simple arithmetics, Eq.(\ref{Section-3-c-48}) can be induced as the same outline of Eq.(\ref{Section-3-c-40}), which completes our proof.

Here provides two methods for determining the concise solution for geodesic distance on tree models $N(t,m)$. While the both computations are proceeded in an iterative manner, the nature concealed by them is completely different from one another. The two techniques, in some sense, have the same impact on calculation process from the complexity point of view, in particular, when an edge serves as seed. If we are in some other situations, for instance, where the seed is assigned as a larger tree $\mathcal{T}(\mathcal{V},\mathcal{E})$ on tens of vertices, the method based on self-similarity seems to become inadequate but our novel ways addressed in this paper can still be adequately employed to work well where only requirement is to know the geodesic distance $\mathcal{S}$ and vertex number $|\mathcal{V}|$ before carrying out our technique. Furthermore, our methods are more effective to implement than some universally studied ones built by matrix, such as, Laplacian spectral and eigenvectors of underlying structure. One of most important reasons for this is the sparsity of adjacency matrix corresponding to tree models of such types. Equivalently, while the total number of entries in adjacency matrix of tree models in question increases exponentially over time, the nonzero ones are order of magnitude as vertex number. Therefore, some matrix methods may be used to calculate the exact solution for geodesic distance on tree models discussed here after executing a number of matrix operations. Sometimes, these such operations appear to be complicated because an arbitrary tree $\mathcal{T}(\mathcal{V},\mathcal{E})$ can be considered as a seed. Facing with situations of such kind, one can choose our technique to derive desired consequences because our methods have no great dependence on the choice of original trees (seeds). To put this another way, they can be quite competent to deal with such tasks in a reasonable time using present laptops when choosing a seed with up to hundreds of vertices.

As before, the average geodesic distance $\langle\mathcal{S}(t,m)\rangle$ on tree model $N(t,m)$ can be readily derived by both definition and Eq.(\ref{Section-3-c-40}), and is exhibited in the below application.

\textbf{Application 4} After $t$ time steps, the solution of average geodesic distance $\langle\mathcal{S}(t,m)\rangle$ on tree model $N(t,m)$ will follow

\begin{equation}\label{Section-3-a-40}
\langle\mathcal{S}(t,m)\rangle\approx \frac{2^{t+1}(m+1)^{2}}{(m+2)(2m+3)}\approx|\mathcal{V}_{t}|^{\gamma}= O(D_{t})
\end{equation}
where exponent $\gamma$ is equal to $\frac{\ln 2}{\ln m+2}$.

Based on statements from Eq.(\ref{Section-3-a-30}) and Eq.(\ref{Section-3-a-40}), we can immediately arrive at the following theorem which says a more general phenomenon.

\textbf{Theorem 4}  Given an arbitrary tree $\mathcal{T}(\mathcal{V},\mathcal{E})$, the exact solution for average geodesic distance $\langle\mathcal{S^{\star}}_{t}\rangle$ on the ($1,m$)-star-fractal tree $\mathcal{T^{\star}}_{t}(\mathcal{V^{\star}}_{t},\mathcal{E^{\star}}_{t})$ is

\begin{equation}\label{Section-3-t-40}
\begin{aligned}\langle\mathcal{S^{\star}}_{t}\rangle&=\frac{2\mathcal{S^{\star}}_{t}}{|\mathcal{V}^{\star}_{t}|(|\mathcal{V}^{\star}_{t}|-1)}\\
&\approx|\mathcal{V}^{\star}_{t}|^{\gamma^{\star}}= O(D^{\star}_{t})
\end{aligned}
\end{equation}
where exponent $\gamma^{\star}=\frac{\ln 2}{\ln m+2}$ in the limit of large graph size, and symbol $D^{\star}_{t}$ is viewed as the diameter of tree $\mathcal{T^{\star}}_{t}(\mathcal{V^{\star}}_{t},\mathcal{E^{\star}}_{t})$. Note that we have used Eqs. (\ref{Section-2-1}) and (\ref{Section-3-c-30}).

In addition, the best studied case of ($1,m$)-star-fractal tree $\mathcal{T^{\star}}_{t}(\mathcal{V^{\star}}_{t},\mathcal{E^{\star}}_{t})$ is in practice the T-graph($t$) which can be created from a single edge as a seed by using ($1,1$)-star-fractal operation until $t$ time step. Technically, the T-graph($t$) may also be reconstructed from three preceding T-graph($t-1$)s ($t\geq1$) by vertex-merging-operation plotted in Fig.5. Apparently, let parameter $m$ in Eq.(\ref{Section-3-c-40}) be equivalent to $1$ and then one may write the next corollary.

\begin{figure}
\centering
  \includegraphics[height=3.3cm]{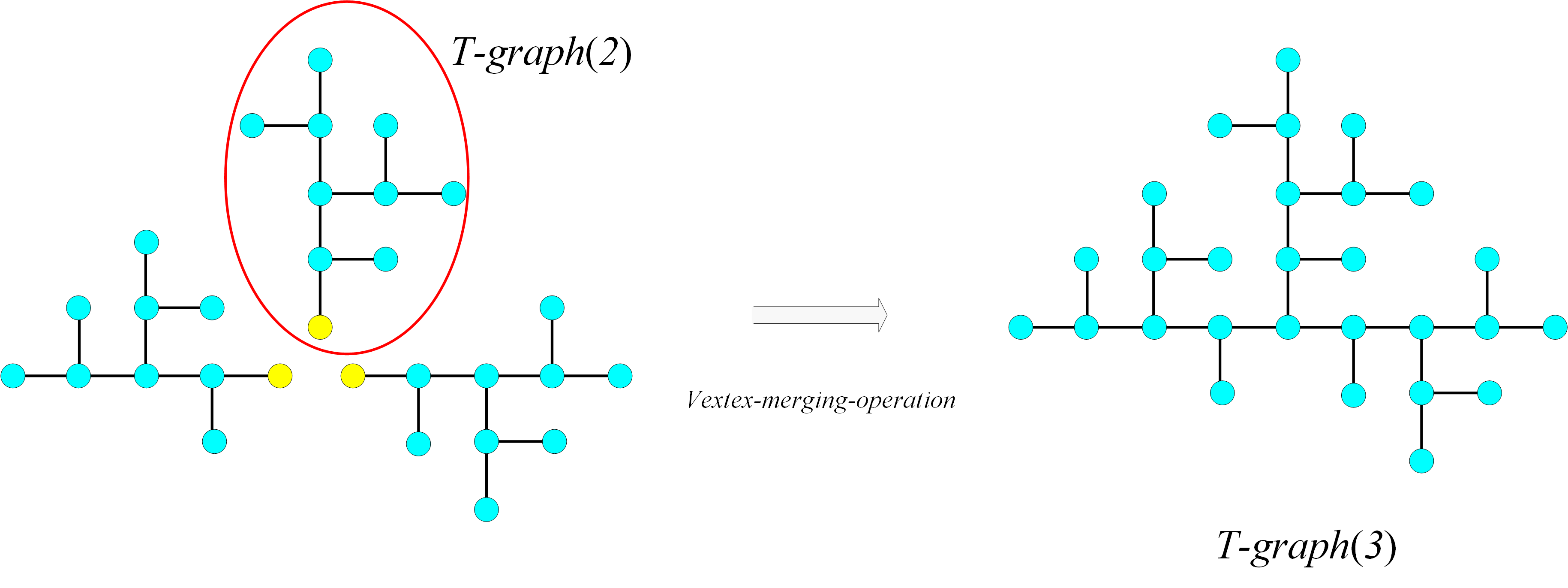}\\
{\small Fig.5. The diagram of T-graph($3$) that is able to be constructed from three T-graph($2$)s by applying vertex-merging-operation among three external vertices (colored in yellow online). An algorithm for producing T-graph($t$) is introduced in Sec.4 of Supplementary Materials.      }
\end{figure}

\emph{Corollary 5} After $t$ time steps, the solution of geodesic distance $\mathcal{S}_{t}$ on the T-graph($t$) will follow
\begin{equation}\label{Section-3-c-50}
\mathcal{S}_{t}=3^{t}+\frac{2^{t+2}+5}{5}\times3^{2t-1}-\frac{3^{t+1}}{5}.
\end{equation}

At the same time, it is not difficult to directly obtain the fifth application.

\textbf{Application 5} After $t$ time steps, the solution of average geodesic distance $\langle\mathcal{S}_{t}\rangle$ on the T-graph($t$) will follow

\begin{equation}\label{Section-3-a-50}
\langle\mathcal{S}_{t}\rangle\approx2^{t}\approx|\mathcal{V}_{t}|^{\gamma^{\star}}, \qquad \gamma^{\star}=\frac{\ln 2}{\ln 3}
\end{equation}
where $|\mathcal{V}_{t}|$ is the total number of vertices of the T-graph($t$).

It is worth noticing that the ($1,m$)-star-fractal trees $\mathcal{T^{\star}}_{t}(\mathcal{V^{\star}}_{t},\mathcal{E^{\star}}_{t})$ all display some similar topological structural properties with each other. These such properties include: (i) Each central vertex added by ($1,m$)-star-fractal operation at any time step $t_{i}$ ($1\leq t_{i}\leq t$) into models $\mathcal{T^{\star}}_{t}(\mathcal{V^{\star}}_{t},\mathcal{E^{\star}}_{t})$ has degree $m+2$, (ii) The total number of leaves of models $\mathcal{T^{\star}}_{t}(\mathcal{V^{\star}}_{t},\mathcal{E^{\star}}_{t})$ keep successively changed over time \cite{Ma-20191}, and (iii) Each member of model family $\mathcal{T^{\star}}_{t}(\mathcal{V^{\star}}_{t},\mathcal{E^{\star}}_{t})$ exhibits homogeneous topological structure. As mentioned above, ($1,m$)-star-fractal trees $\mathcal{T^{\star}}_{t}(\mathcal{V^{\star}}_{t},\mathcal{E^{\star}}_{t})$ just share property (iii) with first-order subdivision trees $\mathcal{T'}_{t}(\mathcal{V'}_{t},\mathcal{E'}_{t})$ and, however, show sharply different appearances from the latter due to another two properties, in particular, property (i). Except for these difference addressed here, there should be other potential ones which will still wait to unveil in the future.

\section{Random walks}

Consider that the above discussions corresponding to average geodesic distance on tree models under consideration are thought of as some evident results by directly applying our methods, then this section will demonstrate a kind of potential applications by means of our methods, i.e., explicitly determining precise solutions of mean first-passage time ($MFPT$) for random walks on the proposed tree models. On the one hand, we indeed obtain desired expressions in perfect agreement with the results previously reported. On the other hand, our methods can perform much better to handle such problems in a general environment than some widely used manners by diminishing redundant computations. To organize the outline of our work narrated below in a self-contained manner, we have to revisit some notations and terminologies relevant to random walks on graph (network).

As the discrete-time representative of Brownian motion and diffusive processes, random walks have proven useful in a wide range of distinct applications in the past \cite{Wu-2016}-\cite{Fouss-2007}. Nonetheless, random walk still keeps quite active at present and hence attracts more attention \cite{Alfred-2018}-\cite{Dragana-2019}. As known, the random walk describes an ideal situation in which a walker (particle) has no information about the underlying graph $\mathcal{G}(\mathcal{V},\mathcal{E})$ and just chooses uniformly at random one vertex of its neighbor set to move on. The random walk of this type, called the unbiased Markov random walks as well, can be depicted by using Markov chains \cite{Kemeny-1976}. For a walker performing random walk on a graph $\mathcal{G}(\mathcal{V},\mathcal{E})$, the most fundamental and significant issue to solve is to analytically determine the first-passage time $FPT$ from source (start vertex) to trap (destination vertex).

In the last several decades, the random walks performed on a graph $\mathcal{G}(\mathcal{V},\mathcal{E})$ with a single trap have been extensively studied \cite{Haynes-2008,Agliari-2008} and hence some interesting results, which reveal some scaling relations and dominating behavior on graph, have been reported. Nevertheless, more and more researchers think that in some real environments each vertex of graph $\mathcal{G}(\mathcal{V},\mathcal{E})$ has very likely to be selected to serve as a trap. It has been proved that the location of traps and underlying topological structure both strongly affect the behavior of random walks. Hence, a quantity called mean first-passage time ($MFPT$) is commonly adopted as a measure for describing the efficiency of random walks with the perfect trap uniformly allocated at all vertices on graph $\mathcal{G}(\mathcal{V},\mathcal{E})$. The smaller the $MFPT$ is, the higher the efficiency is, and vice versa. In view of the obvious importance and ubiquity of random walks themselves, we will study random walks on our models using electrical network $\mathcal{G}^{\otimes}(\mathcal{V}^{\otimes},\mathcal{E}^{\otimes})$ which can be obtained from its underlying graph $\mathcal{G}(\mathcal{V},\mathcal{E})$ by placing a unit resistance on every edge $uv\in \mathcal{E}$.

Given a graph $\mathcal{G}(\mathcal{V},\mathcal{E})$ of interest, the $MFPT$ is by definition equal to the value averaged over first-passage times ($FPT$) of all pairs of vertices $u$ and $v$ in $\mathcal{V}$. Here, the $FPT$ for any vertex pair can be expressed on the basis of the fundamental matrix corresponding to graph $\mathcal{G}(\mathcal{V},\mathcal{E})$, and then the fundamental-matrix method for calculating the $MFPT$ for random walk on graph $\mathcal{G}(\mathcal{V},\mathcal{E})$ is to calculate the inversion of number $|\mathcal{V}|$ of matrices with cardinality $(|\mathcal{V}|-1)\times(|\mathcal{V}|-1)$. This evidently implies that such type of method is just adequately employed for small graphs but becomes prohibitively
difficult to consider some other graphs with thousands of vertices. To address this issue, one can make use of another candidate introduced in \cite{Israel-2003} which is in essence established based on the matrix of graph $\mathcal{G}(\mathcal{V},\mathcal{E})$. Yet, using this method, the issue above to address can easily be induced to calculate the pseudoinverse of the Laplacian matrix \textbf{L} of graph $\mathcal{G}(\mathcal{V},\mathcal{E})$ now, which allows us to compute the $FPT$ between arbitrary pair vertices $u$ and $v$ directly from the inversion of a single $|\mathcal{V}|\times|\mathcal{V}|$ matrix \cite{Rao-1071}. The entry $l_{uv}$ of Laplacian matrix \textbf{L} follows

\begin{equation}\label{Section-4-1}
l_{uv}=
\left\{
\begin{array}{ll}
-1, & \text{an edge connects vertex $u$ to $v$ }\\
\;k_{u}, & \text{$u=v$}\\
\;\;0, & \text{otherwise}\\
\end{array}
\right.
\end{equation}
In the language of matrix theory, the Laplacian matrix \textbf{L} can be compacted as $\mathbf{L}=\mathbf{Z}-\mathbf{A}$  where symbol \textbf{A} is the adjacency matrix of corresponding graph $\mathcal{G}(\mathcal{V},\mathcal{E})$ and $\mathbf{Z}$ is the diagonal matrix that may be defined as follows: the $i$-th diagonal entry is $k_{i}$, while all non-diagonal elements are zero, i.e., $\mathbf{Z}=diag(k_{1}, k_{2}, ..., k_{|\mathcal{V}|})$. And then, the pseudoinverse of the Laplacian matrix \textbf{L} is

\begin{equation}\label{Section-4-2}
\mathbf{L}^{\ast}=\left(\mathbf{L}-\frac{\mathbf{E}\mathbf{E}^{\top}}{|\mathcal{V}|}\right)^{-1}+\frac{\mathbf{E}\mathbf{E}^{\top}}{|\mathcal{V}|}
\end{equation}
where vector $\mathbf{E}=(1,1,1,...,1_{|\mathcal{V}|})^{\top}$. If let $FPT_{uv}$ denote as the first time took by a walker on graph $\mathcal{G}(\mathcal{V},\mathcal{E})$ to arrive at vertex $v$ from its source vertex $u$, then one is able to write

\begin{equation}\label{Section-4-3}
FPT_{uv}=\sum_{i=1}^{|\mathcal{V}|}(l^{\ast}_{ui}-l^{\ast}_{uv}-l^{\ast}_{vi}+l^{\ast}_{vv})l_{ii}
\end{equation}
here $l^{\ast}_{ij}$ is entry of the matrix $\mathbf{L}^{\ast}$ well defined above and $l_{ii}$ the \emph{i}-th entry of the diagonal of the Laplacian matrix $\textbf{L}$. Based on this, the mean first-passage time ($MFPT$) may be expressed as

\begin{equation}\label{Section-4-4}
MFPT=\frac{1}{|\mathcal{V}|(|\mathcal{V}|-1)}\sum_{u\neq v\in \mathcal{V}}\sum_{v=1}^{|\mathcal{V}|}FPT_{uv}.
\end{equation}

By far, Eqs.(\ref{Section-4-3}) and (\ref{Section-4-4}) say that answering $MFPT$ may be reduced to calculate
the entries of the pseudoinverse matrix $\mathbf{L}^{\ast}$ and since its complexity becomes much lighter than the previous case because of only requirement for inverting a $|\mathcal{V}|\times|\mathcal{V}|$ matrix. Nevertheless, as shown above, the total number of vertices of each model studied in this paper increases exponentially over time step $t$, which leads this matrix-based technique to becoming too tedious to obtain an exact formula for $MFPT$. Fortunately, the topological structure of each model and the relationship of effective resistance to $FPT$ together help us to analytically accomplish the computation of $MFPT$ \cite{Tetali-1991}. Specifically, the effective resistance $R_{uv}$ between a pair of vertices $u$ and $v$ in an electrical network $\mathcal{G}^{\otimes}(\mathcal{V}^{\otimes},\mathcal{E}^{\otimes})$ can without difficulty be transformed to calculate
the $FPT_{uv}$ on its corresponding underlying graph $\mathcal{G}(\mathcal{V},\mathcal{E})$, i.e., $R_{uv}=(FPT_{uv}+FPT_{vu})/2|\mathcal{E}|$. In practice, the expression of numerator, $FPT_{uv}+FPT_{vu}$, is customarily viewed as the commute time $C_{uv}$ between vertices $u$ and $v$, that is, $C_{uv}=FPT_{uv}+FPT_{vu}$. In other words, $R_{uv}$ is denoted by $C_{uv}/2|\mathcal{E}|$.

Armed with the statements above, Eq.(\ref{Section-4-4}) may be reorganized as

\begin{equation}\label{Section-4-5}
MFPT=\frac{1}{|\mathcal{V}|}\sum_{u\neq v\in \mathcal{V}}\sum_{v=1}^{|\mathcal{V}|}R_{uv}
\end{equation}
here we already take equality $C_{uv}=C_{vu}=2|\mathcal{E}|R_{uv}$ for any couple of vertices $u$ and $v$. Given a general graph $\mathcal{G}(\mathcal{V},\mathcal{E})$, the complexity of effective resistance computation of its corresponding electrical network $\mathcal{G}^{\otimes}(\mathcal{V}^{\otimes},\mathcal{E}^{\otimes})$ is still to invert a $|\mathcal{V}|\times|\mathcal{V}|$ matrix which is in some extent the same as that said in Eq.(\ref{Section-4-4}). This does not appear to lighten our workload. However, we want to stress that Eq.(\ref{Section-4-5}) is versatile for any graph $\mathcal{G}(\mathcal{V},\mathcal{E})$ of interest. For some graphs with specific topological structure, this transformation addressed in Eq.(\ref{Section-4-5}) is found to work well. For instance, tree highlights the convenience and significance of transformation of such type by supporting a fact that the effective resistance of arbitrary two distinct vertices is completely equivalent to the geodesic distance between this vertex pair. Since then, we are allowed to use the lights shed by Eq.(\ref{Section-4-5}) to find out rigorous expression of $MFPT$ for each member of our models.

From now on, let us divert insights into discussing the mean first-passage time $(MFPT)$ on two types of tree models, $\mathcal{T'}_{t}(\mathcal{V'}_{t},\mathcal{E'}_{t})$ and $\mathcal{T^{\star}}_{t}(\mathcal{V^{\star}}_{t},\mathcal{E^{\star}}_{t})$, using Eq.(\ref{Section-4-5}). The closed-form solutions to $MFPT$ for random walk on them are reported in the next theorems.

\textbf{Theorem 5}  Given an arbitrary tree $\mathcal{T}(\mathcal{V},\mathcal{E})$, the solution for the mean first-passage time $MFPT'_{t}$ on its first-order subdivision tree $\mathcal{T'}_{t}(\mathcal{V'}_{t},\mathcal{E'}_{t})$ will follow

\begin{equation}\label{Section-4-6}
MFPT'_{t} \approx\frac{2^{2t+2}\mathcal{S}}{(|\mathcal{V}|-1)}-\frac{2^{2t+2}}{3}+(2-2^{2t+1})(|\mathcal{V}|-1).
\end{equation}

For a tree $\mathcal{T}(\mathcal{V},\mathcal{E})$ with finite number of vertices, according to Eq.(\ref{Section-3-a-20}), Eq.(\ref{Section-4-6}) may be asymptotically expressed as

\begin{equation}\label{Section-4-7}
MFPT'_{t} \approx2^{2t}|\mathcal{V}|= O(|\mathcal{V'}_{t}|^{\lambda'})
\end{equation}
where exponent $\lambda'$ is equal to $2$ in the limit of large graph size. This means that for first-order subdivision tree $\mathcal{T'}_{t}(\mathcal{V'}_{t},\mathcal{E'}_{t})$, the $MFPT'_{t}$ is a power function with respect to its order. The more the vertex number, the more the $MFPT'_{t}$.  At the moment, it is apparent to find out an equality $\lambda'=1+\gamma'$. By analogous computation to Eq.(\ref{Section-4-6}), one can have the next theorem.

\textbf{Theorem 6}  Given an arbitrary tree $\mathcal{T}(\mathcal{V},\mathcal{E})$, the solution for the mean first-passage time $MFPT^{\star}_{t}$ on its ($1,m$)-star-fractal tree $\mathcal{T^{\star}}_{t}(\mathcal{V^{\star}}_{t},\mathcal{E^{\star}}_{t})$ will follow

\begin{equation}\label{Section-4-8}
\begin{aligned}MFPT^{\star}_{t} &\approx\frac{4(4+2m)^{t}\mathcal{S}}{(|\mathcal{V}|-1)}-\frac{4(4+2m)^{t}(m+1)}{2m+3}\\
&+(2-2^{t+1})(2+m)^{t}(|\mathcal{V}|-1).
\end{aligned}
\end{equation}

As before, for the large value of time $t$ and finite-size tree $\mathcal{T}(\mathcal{V},\mathcal{E})$, the $MFPT^{\star}_{t}$ will have an asymptotical relationship with its order $|\mathcal{V^{\star}}_{t}|$, as follows

\begin{equation}\label{Section-4-9}
MFPT^{\star}_{t} \approx(4+2m)^{t}|\mathcal{V}|= O(|\mathcal{V^{\star}}_{t}|^{\lambda^{\star}})
\end{equation}
where exponent $\lambda^{\star}=1+\gamma^{\star}$ in the large graph size limit. This shows that the $MFPT^{\star}_{t}$ for random walk on ($1,m$)-star-fractal tree $\mathcal{T^{\star}}_{t}(\mathcal{V^{\star}}_{t},\mathcal{E^{\star}}_{t})$ grows as a power-law function of its order $|\mathcal{V^{\star}}_{t}|$ which is in perfect agreement with the result in \cite{Yuan-2010}.
Now let us show a experimental simulation for purpose of detailed comparison between our closed-form solution derived in Theorem 6 and existing results.

\begin{figure}
\centering
  \includegraphics[height=5.7cm]{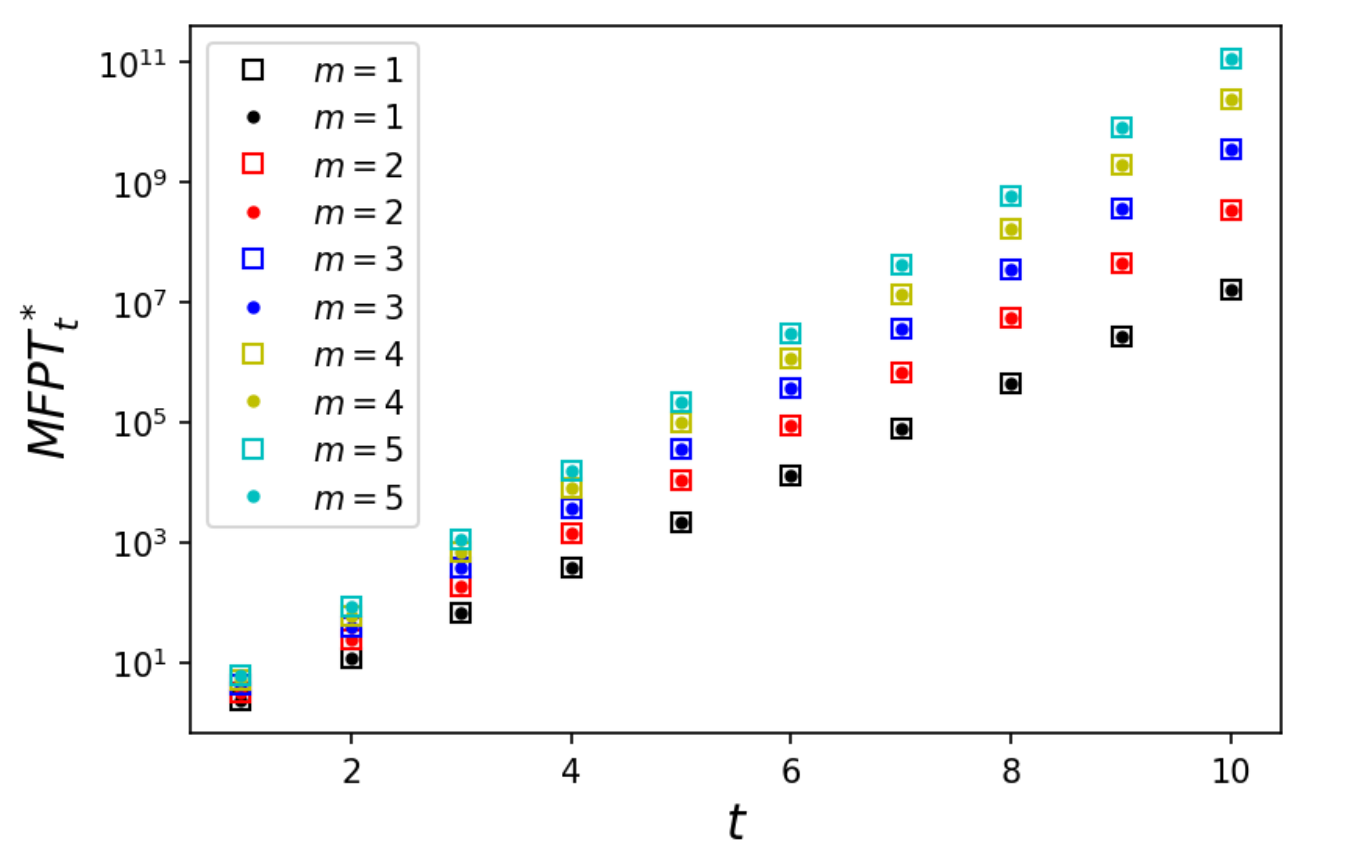}\\
{\small Fig.6. The diagram of comparisons between our derived solutions and previous results. Here, the open squares represent the results obtained using methods reported in \cite{Yuan-2010}, and our consequences are all indicated by solid cycles.    }
\end{figure}

Clearly, as plotted in Fig.6, the results reveal that the exact solutions built by us are in perfect agreement with experimental simulations, implying that our consequences are consistent. On the other hand, the development of our closed-form solutions are more easy to proceed than pre-existing ones including that expressions published in \cite{Yuan-2010}. The most important reason is because there are no more complicated mathematical calculations, for instance, determining all the nonnegative eigenvalues of Laplacian matrix corresponding to model under consideration, occurring in our calculations. In order to clarify the convenience of our methods, interested readers are encouraged to read more details about techniques based on Laplacian matrix for addressing the problems of this type, such as Ref.\cite{Yuan-2010}.

To distinguish how topological structures of the both classes of tree models affect mean first-passage time for random on them, we here need to review some prominent topological structural properties of the both. As known, the special case of the ($1,m$)-star-fractal tree $\mathcal{T^{\star}}_{t}(\mathcal{V^{\star}}_{t},\mathcal{E^{\star}}_{t})$, T-graph, is a fractal with the fractal dimension $d_{f}=\frac{\ln3}{\ln2}$ and the random-walk dimension $d_{w}=\frac{\ln6}{\ln2}=1+d_{f}$. In the meantime, the spectral dimension of T-graph is $\widetilde{d}=\frac{2d_{f}}{d_{w}}=\frac{\ln 9}{\ln 6} <2$, suggesting which a random walk on it is persistent \cite{Erik-2005}. Similarly, our model $\mathcal{T^{\star}}_{t}(\mathcal{V^{\star}}_{t},\mathcal{E^{\star}}_{t})$ has the fractal dimension $d^{\star}_{f}=\frac{\ln(m+2)}{\ln2}$, the random-walk dimension $\lambda^{\star}=\frac{\ln2(m+2)}{\ln2}=1+d_{f}$ and the spectral dimension $\widetilde{d}^{\star}=\frac{2d^{\star}_{f}}{\lambda^{\star}}=\frac{\ln (m+2)^{2}}{\ln 2(m+2)} <2$. Therefore, random walk on our model $\mathcal{T^{\star}}_{t}(\mathcal{V^{\star}}_{t},\mathcal{E^{\star}}_{t})$ is persistent also. By contrast, tree models $\mathcal{T'}_{t}(\mathcal{V'}_{t},\mathcal{E'}_{t})$ do not show fractal phenomena and so have a lack of such rich properties.

\section{Related work}

As mentioned previously, there is a long history of researches on tree not only because tree itself is the most fundamental and simplest complete graph \cite{Bondy-2008} but also since a great number of applications are closely related to structure of tree \cite{Julia-2019},\cite{Peng-2014},\cite{Zhangj-2014}. In particular, a portion of trees with intriguing properties, for instance, fractal feature, have drawn an increasing attention in various scientific fields \cite{Peng-2014},\cite{Yuan-2010},\cite{Zhang-2011}. One of most significant reasons for this is that some real-world networks such as dendrimers and regular hyperbranched polymers can be interpreted as tree models \cite{Wu-2012} in order to (1) study effect from underlying structure on dynamics taking place on them from theoretical point of view \cite{Alfred-2018}-\cite{Takeharu-2018} and (2) mine useful knowledge about information diffusion on them from practical aspect \cite{Wu-2016}-\cite{Fouss-2007}. More generally, there are two topological parameters, geodesic distance and mean first-passage time, commonly used in the process of addressing the issues mentioned above \cite{Nicolas-2012}-\cite{Beveridge-2016}.

In the past, the most used tree models are almost created based on a single edge as seed in an iterative manner like that employed to produce our tree models in this paper. Some of them are in fact constructed via ($1,m$)-star-fractal operation \cite{Yuan-2010},\cite{Agliari-2008}. Using some typical methods including spectral techniques and discrete Green's functions, the both parameters, geodesic distance and mean first-passage time, have been analytically derived \cite{Peng-2014},\cite{Yuan-2010},\cite{Zhang-2011},\cite{Agliari-2008}. Indeed, these works have given helpful guides for better understanding structural properties behind the generated tree models with respect to practical aspects. Nonetheless, those methods previously used are slightly complicated mainly because they are suitable to general graphs including tree. Hence, according to nature of tree structure, we design novel techniques for precisely calculating the closed-form solutions to the above two parameters. Our results derived are perfectly consistent with the published ones. At the same time, our techniques are developed in a combinatorial manner and hence do not involve more complicated arithmetic such as determining the spectrum of the normalized Laplacian matrix corresponding to tree models, implying that the methods built here are easy to understand and convenient to manipulate. To make further progress, the seed used to create our tree models is not limited to a single edge but an arbitrary tree. It should be noted that our techniques are still competent in analytically obtaining the desired results in such a more general situation as reported in sections 3 and 4. Roughly speaking, our work aims at not only addressing some well-studied problems by simplifying computing formulas but also enriching fundamental understanding by generalizing construction of models.

It should be mentioned that there are a large number applications related to tree models proposed here. Among which, many problems in physics and chemistry are connected to random walks on fractal structures such as the special case of our ($1,m$)-star-fractal trees, i.e., T-graph($t$) \cite{Agliari-2008},\cite{Wu-2012}. As a concrete example, the well-known T-graph($t$), as a classic presentative of networks with low-dimensional and fractal structures, has been used to mimic underlying structures of some regular hyperbranched polymers and then to study the scaling law of information diffusion on these polymers via calculating Laplacian spectra of the corresponding Laplacian matrix \cite{Agliari-2008}. In other words, the catch in the problem above is to precisely determine all non-zero eigenvalues of Laplacian matrix. Such a manner using Laplacian spectra has been used in \cite{Wu-2012} in order to address this issue. In fact, this is a common tool for addressing the problem above and thus there are no doubts that some unnecessarily complicated computations might be brought. The computational complexity based on computing non-zero eigenvalues of matrix may be firmly proved in \cite{Chung-1994}. Also, see \cite{Yuan-2010} for more details. As known, the underlying structures of polymers of this type are not completely the same as T-graph($t$). So, we generalized T-graph($t$) using $(1,m)$-star-fractal operation based on an arbitrary tree to generate more general models for purpose of better fitting such type of polymers with varying sizes. As a result, those previous methods, for instance, methods from Laplacian spectra theory, suitable for the typical T-graph($t$) and its some simple variants can no longer be quite useful to precisely derive some topological parameters, such as, mean first-passage time mentioned in Theorem 6, to understand how the underlying structure affects the dynamics taking place in such a condition as considered here. On the other hand, as stated above, the methods proposed in this paper can be clearly adequate for answering these issues. That is to say, we first derive the closed-form solution of geodesic distance on model and then obtain the exact expression of mean first-passage time using Eq.(\ref{Section-4-5}). Therefore, the main purpose of this paper is to fulfill such a gap.

Similarly, our another example trees built using ($1,2$)-star-fractal operation, called Peano basin fractal, have been employed to develop Peano river network for investigating the low connections of natural river channels \cite{Bartolo-2016}. In particular, the authors in \cite{Bartolo-2016} are interesting in the robustness of network based on natural river channels. Taking into account irregular structure of the resulting network itself, they conducted a great number of numerical simulations and finally obtained desirable results. In theory, the mean effective resistance\footnote[3]{For a network $\mathcal{G}(\mathcal{V},\mathcal{E})$, its mean effective resistance is defined as $\langle\mathcal{R}_{\mathcal{G}}\rangle=\mathcal{R}_{\mathcal{G}}/|\mathcal{V}|(|\mathcal{V}|-1)$ where symbol $\mathcal{R}_{\mathcal{G}}$ is the Kirchhoff index of network. }, as a useful measure, may be used to estimate the robustness of underlying structure of network. In this kind of networks established upon ($1,2$)-star-fractal operation, the calculation associated with mean effective resistance can still be converted into computation of non-zero eigenvalues of the corresponding Laplacian matrix \cite{Kemeny-1976}. Clearly, in the sense, our methods are more competitive than many other ones as previously.  Most generally, Patterson \emph{et al} had discussed distributed consensus algorithms in fractal networks including our ($1,m$)-star-fractal trees where agents are subject to external disturbances, and characterized the coherence of these networks in terms of an $H_{2}$ norm of the system that captures how closely agents track the consensus value \cite{Patterson-2011}.

\section{Conclusion and discussion }

To summarize, based on two different types of operations, i.e., first-order subdivision and (1,$m$)-star-fractal operation, we generate two classes of tree models whose seed need not be limited to a single edge but is an arbitrary tree $\mathcal{T}(\mathcal{V},\mathcal{E})$. This implies that our work covers the case of a single edge. In addition, we propose a family of novel and useful enumeration methods for determining the exact solution for geodesic distance on each member of tree models under consideration. In comparison with some commonly used methods for such problems, for instance, Laplacian spectral and eigenvalue, our techniques not only provide what we want to calculate, but also reduce a significant amount of computations, i.e., diminishing the demand of redundant time and space memory when considering this type of tree models with thousands of vertices or more. At the same time, thanks to the special topological structure of models $N(t,m)$, we derive exact formulas for its geodesic distance and average geodesic distance using self-similar method, respectively. This is in strong agreement with those results obtained using our techniques addressed in this paper, further suggesting that the convenience and correctness of our methods.

To highlight the potential applications of our methods, we study the mean first-passage time ($MFPT$) for random walks on tree models built here by connection between random walk and electrical network. Indeed, we obtain the precise expression of $MFPT$ on each tree model. While part of these consequences have been captured in another general manner as in \cite{Yuan-2010}, this does not erase our contributions. As above, two of reasons for this are to have a chance to choose an arbitrary tree as a seed to create our desired models and to cut off a heavy number of computation needs.

However, we want to express that our work is only a tip of the iceberg. On the other hand, we believe that the lights shed by our techniques can be helpful to conduct other related research. Meantime, referring to Ockham's razor we would like to point out that it is important and necessary to develop some new and professional methods for answering many special cases of scientific issues which people attempt to seek for a universal tool to address. Along this research pathes of such type, we will continue to study more general models in the next future.

\section*{ACKNOWLEDGMENT}

The research was supported by the National Key Research and Development Plan under grant 2017YFB1200704 and the National Natural Science Foundation of China under grant No. 61662066.

\vskip 0.5cm

\textbf{Fei Ma} is currently working toward the PhD degree in the School of Electronics Engineering and Computer Science,
Peking University, Beijing, CHINA. His current research interests include Graph theory with applications, Random walks, Complex network, and Privacy-preserving. He is a Reviewer of Mathematical Reviews.

\vskip 0.3cm
\textbf{Ping Wang} received his doctorate degree in Computer Science from the University of Massachusetts, USA in 1996. He is currently a professor at Peking University, China. Dr. Wang has authored or co-authored over 80 papers in journals or proceedings such as IEEE TDSC, ACM CCS, IEEE ICWS and IEEE TII, USENIX Security, etc. His research interests include Internet of Things, Distributed Computing, and Information Security. He is an IEEE senior member.

\vskip 0.3cm
\textbf{Xudong Luo} is now studying for the Ph.D degree in the College of Mathematics and Statistics, Northwest Normal University, CHINA. His interests include Infinite dimensional dynamical systems and Synchronization in complex system.


{\footnotesize
\begin{thebibliography}{9}

\setlength{\parskip}{0pt}


\bibitem{Petrushevski-2019} S. Petrushevski, M. Gusev and V. Zdraveski. ``Calculating average shortest path length using Compute Unified Design Architecture (CUDA)". In Proceedings of IEEE 42nd International Convention on Information and Communication Technology, Electronics and Microelectronics (MIPRO). pp. 186-189, 2019.

\bibitem{Yen-2013}  C.C. Yen, M.Y. Yeh and M.S. Chen. ``An efficient approach to updating closeness centrality and average path length in dynamic networks". In Proceedings of IEEE 13th International Conference on Data Mining. pp. 867-876, 2013.

\bibitem{Spyros-2017} S. Angelopoulos, D. Arsenio and C. Durra. ``Infinite linear programming and online searching with turn cost". Theoretical Computer Science. vol. 670, no. 29, pp. 11-22, 2017.

\bibitem{Jia-2018} J.Y. Li, C. Yang, C.J. Fu, Y.C. Gao and H.C. Yang. ``Cooperative epidemics spreading under resource control". Chaos. vol. 28, no. 11, pp. 113116:1-113116:7, 2018.

\bibitem{Wei-2014} W. Huang, S.Y. Chen and W.L. Wang. ``Navigation in spatial networks: A survey". Physica A. vol. 393, no. 1, pp. 132-154, 2014.

\bibitem{Brodka-2011} P. Brodka, P. Stawiak and P. Kazienko. ``Shortest path discovery in the multi-layered social network". In Proceedings of IEEE International Conference on Advances in Social Networks Analysis and Mining. pp. 497-501, 2011.

\bibitem{Michael-2019} L.P. Wackett. ``Microbiology of produced waters: An annotated selection of World Wide Web sites relevant to the topics in environmental microbiology". Environmental Microbiology. vol. 21, no. 1, pp. 1511-1512, 2019.

\bibitem{Michael-2017} M. Golosovsky and S. Solomon. ``Growing complex network of citations of scientific papers: Modeling and measurements". Phys. Rev. E., vol. 95, no.1, pp. 012324:1-012324:9, 2017.

\bibitem{Ulrich-2019} U, Brose, P. Archambault, A.D. Barnes, et al. ``Predator traits determine food-web architecture across ecosystems". Nature Ecology Evolution. vol. 3, no. 6, pp. 919-927, 2019.

\bibitem{Julia-2019} J.P. Casasnovas, M. Gerlach, N. Aguirre and L.A.N. Amaral. ``Large-scale analysis of micro-level citation patterns reveals nuanced selection criteria". Nature Human Behaviour. vol. 3, no. 6, pp. 568-575, 2019.

\bibitem{MaF-2019} F. Ma, P. Wang and B. Yao. ``Random walks on Fibonacci treelike models: emergence of power law". To submit. arXiv:1904.11314v1.

\bibitem{Luo-2019} X.D. Luo, F. Ma and W.T. Xu. ``Exact solutions for geodesic distance on treelike models with some constraints". To submit. arXiv:1909.07041.

\bibitem{Peng-2014} J.H. Peng and G.A. Xu. ``Efficiency analysis of diffusion on T-fractals in the sense of random walks". J. Chem. Phys., vol. 140, no. 13, pp. 134102:1-134102:12, 2014.

\bibitem{Yuan-2010} Y. Lin, B. Wu and Z.Z. Zhang. ``Determining mean first-passage time on a class of treelike regular fractals". Phys. Rev. E., vol. 82, no. 3, pp. 031140:1-031140:12, 2010.

\bibitem{Bondy-2008} J.A. Bondy and U.S.R. Murty. Graph Theory. Springer. 2008.

\bibitem{Zhangj-2014} J. Zhang, E.H. Yang and J.C. Kieffer. ``A universal grammar-based code for lossless compression of binary trees". IEEE Transactions on information theory. vol. 60, no. 3, pp. 1373-1386, 2014.

\bibitem{Zhang-2011} Z.Z. Zhang, B. Wu and G.R. Chen. ``Complete spectrum of the stochastic master equation for random walks on treelike fractals". Europhysics Letters. vol. 96, no. 4, pp. 40009:1-40009:6, 2011.

\bibitem{Aaron-2008} A. Clauset, C. Moore and M.E.J. Newman. ``Hierarchical structure and the prediction of missing links in networks". Nature. vol. 453, no. 7191, pp. 98-101, 2008.

\bibitem{Stefano-2018} S.D. Palma and P.L. Lanzi. ``Traditional wisdom and Monte Carlo tree search face-to-face in the card game scopone". IEEE Transactions on Games. vol. 10, no. 3, pp. 317-332, 2018.

\bibitem{Deng-2019} J. Deng, Q.Q. Ye and Q. Wang. ``Weighted average geodesic distance of Vicsek network". Physica A. vol. 527, no. 8, pp. 121327:1-121327:6, 2019.

\bibitem{Peng-2018} J.H. Peng, G.A. Xu, R.X. Shao, L. Chen and H.E. Stanley. ``Analysis of fluctuations in the first return times of random walks on regular branched networks". J. Chem. Phys., vol. 149, no. 2, pp. 024903:1-024903:6. 2018.

\bibitem{Nobutoshi-2019} N. Ikeda. ``Growth model for fractal scale-free networks generated by a random walk". Physica A. vol. 521, no. 5, pp. 424-434, 2019.

\bibitem{Ma-2019} More generally, the ration $\varpi'_{t}$ of the leaf number $|\mathcal{L'}_{t}|$ of model $\mathcal{T'}_{t}(\mathcal{V'}_{t},\mathcal{E'}_{t})$ and its order $|\mathcal{V'}_{t}|$ will tend to be zero for the large value of time $t$, i.e., $\varpi'_{t}=\lim_{t\rightarrow\infty} \frac{|\mathcal{L'}_{t}|}{|\mathcal{V'}_{t}|}\rightarrow 0$.

\bibitem{Chao-2005} C.M. Song, S. Havlin and H.A. Makse. ``Self-similarity of complex networks". Nature. vol. 433, no. 7024, pp.  392-395, 2005.

\bibitem{Serrano-2008} M.A. Serrano, D. Krioukov and M. Boguna. ``Self-similarity of complex networks and hidden metric spaces". Phys. Rev. Lett., vol. 100, no. 4, pp. 078701:1-078701:4, 2008.

\bibitem{Fei-2018} F. Ma and B. Yao. ``An iteration method for computing the total number of spanning trees and its applications in graph theory". Theoretical Computer Science. vol. 708, no. 1, pp. 46-57, 2018.


\bibitem{Ma-20191} Generally speaking, the ration $\varpi^{\star}_{t}$ of the leaf number $|\mathcal{L^{\star}}_{t}|$ of model $\mathcal{T^{\star}}_{t}(\mathcal{V^{\star}}_{t},\mathcal{E^{\star}}_{t})$ and its order $|\mathcal{V^{\star}}_{t}|$ will tend to be some constant dependent on paremeter $m$ of the ($1,m$)-star-fractal operation for the large value of time $t$, i.e., $\varpi^{\star}_{t}=\lim_{t\rightarrow\infty} \frac{|\mathcal{L^{\star}}_{t}|}{|\mathcal{V^{\star}}_{t}|}\rightarrow \frac{m}{m+1}$. Obviously, the minimal value for $\varpi^{\star}_{t}$ can be arrived at point $m=0$ which is in perfect agreement with that of \cite{Ma-2019}. The upper bound of $\varpi^{\star}_{t}$ is $1$ as $m\rightarrow\infty$. Meantime, there is a relationship close connecting ratio $\varpi^{\star}_{t}$ with exponent $\gamma^{\star}$, namely, $\gamma^{\star}=\frac{\ln 2}{\ln (1+\varpi^{\star}_{t})-\ln \varpi^{\star}_{t}}$, for the large value of parameter $m$.


\bibitem{Wu-2016} Y.B. Wu, R.M. Jin and X. Zhang. ``Efficient and exact local search for random walk based top-k proximity query in large graphs". IEEE Transactions on Knowledge and Data Engineering. vol. 28, no. 5, pp. 1160-1174, 2016.

\bibitem{Janeja-2008} V.P. Janeja and V. Atluri. ``Random walks to identify anomalous free-form spatial scan windows".  IEEE Transactions on Knowledge and Data Engineering. vol. 20, no. 10, pp. 1378-1392, 2008.

\bibitem{Fouss-2007} F. Fouss, A. Pirotte, J.M. Renders and M. Saerens. ``Random-walk computation of similarities between nodes of a graph with application to collaborative recommendation". IEEE Transactions on Knowledge and Data Engineering. vol. 19, no. 3, pp. 355-369, 2007.


\bibitem{Alfred-2018} A.W. Hales. ``Random walks on visible points". IEEE Transactions on information theory. vol. 64, no. 4, pp. 3150-3152, 2018.

\bibitem{Marat-2012} M.V. Burnashev and A. Tchamkerten. ``Estimating a random walk first-passage time from noisy or delayed observations". IEEE Transactions on information theory. vol. 58, no. 7, pp. 4230-4243, 2012.

\bibitem{Takeharu-2018} T. Shiraga, Y. Yamauchi, S. Kijima and M. Yamashita. ``Deterministic random walks for rapidly mixing chains". SIAM Journal on Discrete Mathematics. vol. 32, no. 3, pp. 2180-2193, 2018.


\bibitem{Rushabh-2016} R. Patel, A. Carron and F. Bullo. ``The Hitting Time of Multiple Random Walks". SIAM Journal on Matrix Analysis and Applications. vol. 37, no. 3, pp. 933-954, 2016.

\bibitem{Dragana-2019} D. Bajovi\'{c}, J.M.F. Moura and D. Vukobratovi\'{c}. ``Detecting random walks on graphs with heterogeneous sensors". IEEE Transactions on information theory. vol. 65, no. 8, pp. 4893-4914, 2019.

\bibitem{Kemeny-1976} J.G. Kemeny and J.L. Snell.  Finite Markov Chains. Springer. 1976.

\bibitem{Haynes-2008}  C.P. Haynesand and A.P. Roberts. ``Global first-passage times of fractal lattices". Phys. Rev. E., vol. 78, no. 4, pp. 041111:1-041111:9, 2008.

\bibitem{Agliari-2008} E. Agliari. ``Exact mean first-passage time on the T-graph". Phys. Rev. E., vol. 77, no. 1, pp. 011128:1-011128:6, 2008.

\bibitem{Israel-2003} A. Ben-Israel and T. Greville. Generalized Inverses: Theory and Applications. Springer. 2003.

\bibitem{Rao-1071} C. Rao and S. Mitra. Generalized Inverse of Matrices and its Applications. Wiley. 1971.

\bibitem{Tetali-1991} P. Tetali. ``Random walks and the effective resistance of networks". J. Theor. Probab., vol. 4, no. 1, pp. 101-109, 1991.

\bibitem{Erik-2005} E.M. Bollt and D. Ben-Avraham. ``What is special about diffusion on scale-free nets". New J. Phys., vol. 7, no. 1, pp. 26-46, 2005.

\bibitem{Wu-2012} B. Wu, Y. Lin, Z.Z. Zhang and G.R. Chen. ``Trapping in dendrimers and regular hyperbranched polymers". J. Chem. Phys., vol. 137, no. 4, pp. 044903:1-044903:7, 2012.

\bibitem{Nicolas-2012} B. Nicolas and H. Cecilia. ``The total path length of split trees". The Annals of Applied Probability. vol. 22, no. 5, pp. 1745-1777, 2012.

\bibitem{Athreya-2017} S. Athreya,  W. L\"{o}hr and A. Winter. ``Invariance principle for variable speed random walks on trees". Annals of Probability. vol. 45, no. 2, pp. 625-667, 2017.

\bibitem{Beveridge-2016} A. Beveridge and J. Youngblood. ``The best mixing time for random walks on trees". Graphs and Combinatorics. vol. 32, no. 6, pp. 2211-2239, 2016.

\bibitem{Chung-1994} F.R.K. Chung. ``Spectral graph theory". American Mathematical Society.

\bibitem{Bartolo-2016} S.D. Bartolo, F. Dell'Accio, G. Frandina, G. Moretti, S. Orlandini and M. Veltri. ``Relation between grid, channel, and Peano networks in high-resolution digital elevation models". Water Resources Research. vol. 52, no. 5, pp. 3527-3546, 2016.

\bibitem{Kemeny-1976} J.G. Kemeny and J.L. Snell. ``Finite Markov Chains". New York, NY, USA: Springer, 1976.

\bibitem{Patterson-2011} S. Patterson and B. Bamieh. ``Network coherence in fractal graphs". In Proceedings of 50th IEEE Conference on Decision and Control and European Control Conference (CDC-ECC). pp. 6445-6450, 2011.




\end{thebibliography}
}
\end{document}